\documentclass[12pt]{article}
\usepackage{amsmath}
\begin{document} 

\begin{center}
{\large{\bf POISSON SUMMATION FORMULA\\FOR THE SPACE OF FUNCTIONALS}}
\end{center}

\bigskip

\begin{center}
\textsc{Takashi NITTA and Tomoko OKADA}
\end{center}

\bigskip
\bigskip
                                
\begin{abstract}
In our last work, we formulate a Fourier
transformation on the infinite-dimensional space of functionals. Here we first
calculate the Fourier
transformation of infinite-dimensional Gaussian distribution $\exp\left(-\pi
\xi\int_{-\infty}^{\infty}\alpha ^2(t)dt\right)$ for
$\xi\in{\bf C}$ with Re$(\xi)>0$, $\alpha \in L^2({\bf R})$, using our formulated
Feynman path integral. Secondly we develop the Poisson summation formula for the
space of functionals, and define a functional $Z_s$, $s\in {\bf C}$, the
Feynman path integral of that corresponds to the Riemann zeta function in the
case Re$(s)>1$.
\end{abstract}

\bigskip
\bigskip
                        
{\bf 0. Introduction}

\medskip

Feynman([F-H]) used the concept of his path integral for physical quantizations.
The word $''$physical quantizations$''$ has two meanings : one is for quantum
mechanics and the other is for quantum field theory. We usually use the
same word $''$Feynman path integral$''$. However the meanings
included in $''$Feynman path integral$''$ are two sides, according to the above.
One is of quantum mechanics and the other is of quantum field theory.
The first Feynman path integral corresponds to a study of functional analysis on
the space of functions.  For functional analysis, there
exist many works from standard analysis and nonstandard analysis. However an
approach has been hard from standard analysis or nonstandard analysis to study the
space of
$''$functionals$''$ associating with the second Feynman path integral. 
 
In our last paper([N-O2]), we defined a delta functional $\delta$ and an
infinitesimal Fourier transformation $F$ in the space of functionals as one of
generalizations for Kinoshita's infinitesimal
Fourier transformation in the space of functions. Historically, in 1962, Gaishi
Takeuchi([T]) introduced an infinitesimal $\delta$-function for the space of
functions under nonstandard analysis.  In 1988, 1990, Kinoshita([K1],[K2]) defined
an infinitesimal Fourier transformation for the space of functions. Nitta and
Okada([N-O1],[N-O2]) defined, for funtionals, an infinitesimal Fourier
transformation, using a concept of double infinitesimal, and calculated the
infinitesimal Fourier transformation for two typical examples. The main idea is to
use the concept of double infinitesimals and putting standard parts twice st(st( .
)). In our theory, the infinitesimal Fourier transformation of $\delta$,
$\delta^2$, ... , and
$\sqrt{\delta}$, ... can be calculated as constant functionals, $1$, infinite, ...
, and infinitesimal, ... . 

Now let $H$ be an even infinite number
in $\,^{\ast}{\bf R }$, and $L$ be an
infinitesimal lattice

\medskip
 
$L:=\left\{\varepsilon z \,\left|\, z \in \, 
^{\ast}{\bf Z},\, -\frac H2 \le
\varepsilon z < \frac H2 \right\},\right.$ where
$\varepsilon = \frac 1H$, and let $H'$
be an even infinite number in $\,^{\star}(\,^{\ast}{\bf R
})$, and $L'$ be an
infinitesimal lattice 

\medskip

$L' :=\left\{\varepsilon' z' \,\left |\,z' \in \, ^{\star}(\,^{\ast}{\bf
Z}),\, -\frac {H'}{2} \le \varepsilon ' z' <
\frac {H'}{2} \right\}\right.$, where
$\varepsilon' = \frac {1}{H'}$.

\medskip
\noindent
Then we extend the calculation in our
previous work to the case of

\medskip

$g_{\xi}(a)=\exp\left(-\pi
\xi \,^{\star}{\varepsilon}\sum_{\varepsilon n\in L}a 
(\varepsilon n)^2\right)$ for
$\xi\in{\bf C}$ with Re$(\xi)>0$. 

\medskip
\noindent
If there exists $\alpha \in L^2$ so that $a
=\,^{\star}\alpha$, then st$\left(\exp\left(-\pi
\xi\,^{\star}{\varepsilon}\sum_{\varepsilon n\in L}a (\varepsilon
n)^2)\right)\right)$ is equal to $\exp\left(-\pi
\xi\int_{-\infty}^{\infty}\alpha ^2(t)dt\right)$. The standard part of the
functional\\
$\exp\left(-\pi
\xi\,^{\star}{\varepsilon}\sum_{\varepsilon n\in L}a 
(\varepsilon n)^2\right)$
corresponds naturally to $\exp\left(-\pi
\xi\int_{-\infty}^{\infty}\alpha ^2(t)dt\right)$ in the standard
meaning. Our Fourier transformation of $\exp\left(-\pi
\xi\,^{\star}{\varepsilon}\sum_{\varepsilon n\in L}a
(\varepsilon n)^2\right)$ is $C_{\xi}(b)\exp\left(-\pi
\xi^{-1} \,^{\star}{\varepsilon}\sum_{\varepsilon n\in L}b 
(\varepsilon n)^2\right)$
where $\textrm{st}(C_{\xi}(b))=\left(\ast\left(
\frac{1}{\sqrt{\xi}}\right)\right)^{H^2}\;(\in
\,^{\ast}{\bf R})$. In the calculation, we assume that the real part $\xi$
is positive. Even if
$\xi$ is
$i$ or
$-i$, the coefficient is equal to
$\left(\ast\left(\frac{1}{\sqrt{\xi}}\right)\right)^{H^2}$ shown in the
previous paper. Furthermore, let $m$ be an integer so that
$m|2\,^{\star}HH'{}^2$, and

\medskip

$g_{im}(a)=\exp(-i\pi m\,^{\star}\varepsilon\sum_{k\in L}a^2(k))$ 

\medskip
\noindent
that associates with $\exp\left(-im\pi
\int_{-\infty}^{\infty}\alpha ^2(t)dt\right)$ in the standard
meaning. Then
we calculate our Fourier transformation for
$g_{im}$ : $(Fg_{im})(b)=C_{im}(b)g_{\frac{1}{im}}(b)$.

\medskip
\noindent
We show that 
$C_{im}(b)=\left(\sqrt{\frac
m2}\frac{1+i^{\frac{2\,^{\star}HH'{}^2}{m}}}{1+i}\right)^
{(\,^{\star}H)^2}$ for positive $m$ and $C_{im}(b)=\left(\sqrt{\frac
{-m}{2}}\frac{1+(-i)^{\frac{2\,^{\star}HH'{}^2}{-m}}}{1-i}\right)
^{(\,^{\star}H)^2}$
for negative $m$ if $m|\frac{b(k)}{\varepsilon'}$ for arbitrary $k$ in
$L$. 

Furthermore Using the second infinitesimal and the
lattice, we extend the Poisson summation formula of finite group to 
infinitesimal Fourier transformations for the space of functions and also
for the space of functionals. For an example, we apply the Poisson
summation formula to the above functional
$g_{\xi}$. If the groups are special, it appears
the
$H^2$-th product of $\theta$-functions, or the constant 
$\left(\ast\left(\frac{1}{\sqrt{\xi}}\right)\right)^{H^2}$. We also apply it to
the functional $g_{im}$. Finally we define a functional that associates to the
Riemann zeta function. Using our Poisson summation formula for functionals, we
study a relationship between the functional and the
Riemann zeta function.

\bigskip
                        
{\bf 1. Preliminaries}

\medskip
 
{\bf 1-1. Infinitesimal Fourier transformations by Kinoshita} (cf.
[Ki],[N-O1],[N-O2])

\medskip

Let $\Lambda$ be an infinite set. Let $F$ be an ultrafilter on $\Lambda$. 
For each $\lambda\in\Lambda$, let $S_{\lambda}$ be a set. We put an equivalence
relation $\sim$ induced from $F$ on $\prod_{\lambda \in\Lambda}S_{\lambda}.$ For
$\alpha =(\alpha_{\lambda}),\; \beta =(\beta_{\lambda})\;(\lambda\in\Lambda )$,

\medskip

$\alpha \sim \beta  \Longleftrightarrow  \{\lambda\in\Lambda
\,|\,\alpha_{\lambda}=\beta_{\lambda}\} \in F$. 

\medskip
\noindent
The set of equivalence classes is called {\it ultraproduct} of $S_{\lambda}$ for $F$
with respect to $\sim$. If $S_{\lambda}=S$ for $\lambda\in\Lambda $, then it is
called {\it ultraproduct} of $S$ for $F$ and it is written as $\,^{\ast}S$. The set
$S$ is naturally embedded in $\,^{\ast}S$ by the following mapping :

\medskip
$s\; (\in S)\mapsto [(s_{\lambda}=s),\lambda\in\Lambda]\; 
(\in \,^{\ast}S)$,

\medskip
\noindent
where
$[\;\;\;]$ denotes the equivalence class with respect to the ultrafilter $F$. We write
the mapping as $\ast$, and call it naturally elementary embedding. From now on, we
identify the image $\ast (S)$ as $S$. 

\medskip

Let $H\; (\in \,^{\ast}{\bf Z})$ be an infinite even number. The infinite number $H$
is even, when for $H=[(H_{\lambda}),\,\lambda\in\Lambda
]$, $\{\lambda\in\Lambda\,|\,H_{\lambda}\textrm{ is even}\}\in F$. We denote
$\frac 1H$ by $\varepsilon$. We define an infinitesimal lattice space
${\bf L}$, an infinitesimal lattice subspace $L$ and a space of functions $R(L)$ on
$L$ as follows :

\medskip

${\bf L}:=\varepsilon\,^{\ast}{\bf Z}=\{\varepsilon z \,|\, z \in \, ^{\ast}{\bf
Z}\}$, 

$L:=\left\{\varepsilon z \,\left|\, z \in \, ^{\ast}{\bf Z},\, -\frac H2 \le
\varepsilon z < \frac H2 \right\}(\subset {\bf L}),\right.$

$R(L):=\left\{\varphi \,\left|\, \varphi \textrm{ is an internal function from } L
\textrm{ to }^{\ast}{\bf C}\right\}.\right.$

\medskip
\noindent
We extend $R(L)$ to the space of periodic functions on ${\bf L}$ with period $H$. We write the same notation $R(L)$ for the space of periodic
functions.

Gaishi Takeuchi([T]) introduced an infinitesimal $\delta$ function. Furthermore
Moto-o Kinoshita ([Ki]) constructed an infinitesimal Fourier transformation
theory on $R(L)$. 

We explain it briefly.
For $\varphi ,\, \psi$ $\in$ $R(L)$, the infinitesimal $\delta$ function, 
the infinitesimal Fourier transformation $F\varphi\;(\in R(L))$, the
inverse infinitesimal Fourier transformation $\overline F\varphi\;(\in
R(L))$ and the convolution
$\varphi\ast \psi\;(\in R(L))$ are defined as follows :

\begin{equation*}
\delta\in R(L),\;\;
\delta (x):= 
\begin{cases}
H & \text{$(x=0)$},\\
0  & \text{$(x \ne 0)$},
\end{cases}
\end{equation*}

$(F\varphi)(p) := \sum _{x \in L} \varepsilon  \exp \left(-2\pi i
px\right)\varphi(x)$,

$(\overline{ F}\varphi)(p) := \sum _{x \in L}
\varepsilon   \exp \left(2\pi i px\right)\varphi(x)$, 

$(\varphi\ast \psi)(x) :=
\sum _{y \in L} \varepsilon \varphi(x-y)\psi(y)$.

\bigskip                     
 
{\bf 1-2. Formulation of infinitesimal Fourier transformation on the space of 
functionals} (cf. [N-O1],[N-O2])

\medskip

To treat a $\,^{\ast}$-unbounded functional $f$ in the nonstandard analysis, we
need a second nonstandardization. Let
$F_2:=F$ be a nonprincipal ultrafilter on an infinite set $\Lambda_2:=\Lambda$ as
above. Denote the ultraproduct of a set $S$ with respect to $F_2$ by $\,^{\ast}S$
as above. Let $F_1$ be another nonprincipal ultrafilter on an infinite set
$\Lambda_1$. Take the $\,^{\ast}$-ultrafilter $\,^{\ast}F_1$ on
$\,^{\ast}\Lambda_1$. For an internal set $S$ in the sense of
$\,^{\ast}$-nonstandardization, let $\,^{\star}S$ be the $\,^{\ast}$-ultraproduct
of $S$ with respect to $\,^{\ast}F_1$. Thus, we define a double ultraproduct
$^{\star}(^{\ast}{\bf R})$, $^{\star}(^{\ast}{\bf Z})$, etc for the set ${\bf
R}$, ${\bf Z}$, etc. It is shown easily that

$$^{\star}(^{\ast}{\bf S})=S^{\Lambda_1\times\Lambda_2}/F_1^{F_2},$$

\medskip
\noindent
where $F_1^{F_2}$ denotes the ultrafilter on $\Lambda_1\times\Lambda_2$ such
that for any $A\subset\Lambda_1\times\Lambda_2$, $A\in F_1^{F_2}$ if and only if 

$$\{\lambda\in\Lambda_1\,|\,\{\mu\in\Lambda_2\,|\,(\lambda,\mu)\in A\}\in
F_2\}\in F_1.$$

\medskip
\noindent
We always work with this double nonstandardization. The natural imbedding
$\,^{\star}S$ of an internal element $S$ which is not considered as a set in
$\,^{\ast}$-nonstandardization is often denoted simply by $S$.

An infinite number in $^{\star}(^{\ast}{\bf R})$ is defined to be greater
than any element in $^{\ast}{\bf R}$. We remark that an infinite number in 
$\,^{\ast}{\bf R}$ is not infinite in $\,^{\star}(\,^{\ast}{\bf R})$, that
is, the word $''$an infinite number in $^{\star}(^{\ast}{\bf R})''$ has a double
meaning. An infinitesimal number in
$^{\star}(^{\ast} {\bf R})$ is also defined to be nonzero and whose absolute
value is less than each positive number in $^{\ast}{\bf R}$.

\medskip   

\textsc{Definition} 1.1.\;\; Let $H (\in \,^{\ast}{\bf Z}),\, H' (\in
\,^{\star}(^{\ast}{\bf Z}))$ be even positive numbers such that $H'$ is larger
than any element in $\,^{\ast}{\bf Z}$, and let $\varepsilon (\in
\,^{\ast}{\bf R}),\,\varepsilon ' (\in \,^{\star}(^{\ast}{\bf R}))$ be 
infinitesimals satifying $\varepsilon H =1,\,\varepsilon ' H' =1$. We define as
follows :

\medskip 

${\bf L} :=\varepsilon \, ^{\ast}{\bf Z}= \{\varepsilon z \,|\,z \in \, 
^{\ast}{\bf
Z}\},\;\;
{\bf L}' :=\varepsilon ' \, ^{\star}(\,^{\ast}{\bf Z})=\{\varepsilon ' z'
\,|\, z' \in \, ^{\star}(\,^{\ast}{\bf Z})\},$
 
$L := \left\{\varepsilon z \,
\left|\,z \in \, ^{\ast}{\bf Z},\, -\frac H2 \le
\varepsilon z < \frac H2 \right\}\right.(\subset {\bf L})$,

$L' :=\left\{\varepsilon' z' \,\left |\,z' \in \, ^{\star}(\,^{\ast}{\bf Z}),\,
-\frac {H'}{2} \le \varepsilon ' z' <
\frac {H'}{2} \right\}(\subset {\bf L}').\right.$

\medskip 
\noindent
We define a latticed space of functions $X$ as follows,

\medskip 

$X := \{a \,|\, a \textrm{ is an internal function with a double meaning,
from}
\star(L)\textrm{ to } L' \}$

\medskip 
\noindent
We define three
equivalence relations $\sim_H$, $\sim_{\star(H)}$ and
$\sim_{H'}$ on
${\bf L}$, $\star({\bf L})$ and ${\bf L}'$ :

\medskip 

$x\sim_H y\Longleftrightarrow x-y\in H\,^{\ast}{\bf Z},\;\;
x\sim_{\star(H)} y\Longleftrightarrow x-y\in \star(H)\,^{\star}(\,^{\ast}{\bf Z}),$

$x\sim_{H'} y\Longleftrightarrow x-y\in H'\,^{\star}(\,^{\ast}{\bf Z})$.

\medskip 
\noindent
Then we identify ${\bf L}/\sim_H$, $\star({\bf L})/\sim_{\star(H)}$ and ${\bf
L}'/\sim_{H'}$ as $L$, $\star(L)$ and $L'$. Since $\star(L)$ is identified 
with $L$, the set 
$\star({\bf L})/\sim_{\star(H)}$ is identified with ${\bf L}/\sim_H$. Furthermore
we represent
$X$ as the following internal set : 

\medskip 
\noindent
$\{a\,|\,a\textrm{ is an internal function with a double meaning,
from }\star({\bf L})/\sim_{\star(H)}\textrm{to }{\bf
L}'/\sim_{H'}\}.$

\medskip 
\noindent
We use the same notation as a function from $\star(L)$
to $L'$ to represent a function in the above internal set. We define the space $A$
of functionals as follows :

\medskip 
$A:= \{ f \,|\, f \textrm{ is an internal function with a double meaning, from
} X 
\textrm{ to } \,^{\star}(\,^{\ast}{\bf C})
\}.$

\medskip

We define an infinitesimal delta function $\delta (a)(\in
A)$, an infinitesimal Fourier transformation of
$f(\in A)$, an inverse infinitesimal Fourier transformation of $f$ and a
convolution of $f$,
$g(\in A)$, by the following :

\medskip

\textsc{Definition} 1.2.
\begin{equation*}
\delta (a):= 
\begin{cases}
(H')^{(\,^{\star}H)^2} & \text{$(a=0)$}, \\
0  & \text{$(a \ne 0)$},
\end{cases}
\end{equation*}

$\varepsilon _0 := (H')^{-{(\,^{\star}H)^2}}
\in \,^{\star}(^{\ast}{\bf R})$,

$(Ff)(b) := \sum _{a \in X} \varepsilon _0  \exp \left(-2\pi i \sum
_{k\in L} a(k)b(k)\right)f(a)$,

$(\overline{ F} f)(b) := \sum _{a \in X}
\varepsilon _0  \exp \left(2\pi i \sum _{k\in L} a(k)b(k)\right)f(a)$,

$(f\ast g)(a)
:=\sum _{a' \in X} \varepsilon _0  f(a-a')g(a')$.

\medskip 
\noindent
We define an inner product on $A$ : 

\medskip 

$(f,g):=\sum_{b\in
X}\varepsilon_0\overline{f(b)}g(b)$, where $\overline{f(b)}$ is the complex conjugate
of $f(b)$. 

\medskip

Replacing the definitions of $L'$, $\delta$, $\varepsilon_0$, $F$,
$\overline{F}$ in  Definition 1.1 and Definition 1.2 by the following, we
shall define another type of infinitesimal Fourier transformation. The
different point is only the definition of an inner product of the space
of functions $X$. In Definition 1.2, the
inner product of $a,\,b(\in X)$ is $\sum _{k\in L}a(k)b(k)$, and in the
following definition, it is 
$\,^{\star}\varepsilon\sum _{k\in L}a(k)b(k)$.

\medskip

\textsc{Definition} 1.3.

\medskip

$L' := \left\{\varepsilon ' z' \,\left|\, z'\in \, ^{\star}(\,^{\ast}{\bf Z}) ,\, -\,
^{\star}H \frac {H'}{2} \le \varepsilon ' z' < \, ^{\star}H \frac {H'}{2}
\right\},\right.$
\begin{equation*}
\delta (a):= 
\begin{cases}
(\,^{\star}H)^{\frac{(\,^{\star}H)^2}{2}}H'{}^{(\,^{\star}H)^2} &
\text{$(a=0)$}, \\ 0  & \text{$(a \ne 0)$},
\end{cases}
\end{equation*}

$\varepsilon_0
:=(\,^{\star}H)^{-\frac{(\,^{\star}H)^2}{2}}H'{}^{-(\,^{\star}H)^2}$

$(Ff)(b) := \sum _{a \in X} \varepsilon _0  \exp
\left(-2\pi i \,^{\star}\!\varepsilon \sum _{k\in L} a(k)b(k)\right)f(a)$,

$(\overline{F} f)(b) :=\sum _{a \in X} \varepsilon _0  \exp \left(2\pi i
\,^{\star}\!\varepsilon \sum _{k\in L} a(k)b(k)\right)f(a).$

\medskip

Then we obtain the following theorem :

\medskip

\textsc{Theorem} 1.4([N-O2]).

$\textrm{(1)}\;\; \delta = F1=\overline{F} 1,\;\;
\textrm{(2)}\;\; F \: \textrm{is unitary}, \, F^4 =1, \overline{ F} F =F
\overline{ F}=
1$,

$\textrm{(3)}\;\; f\ast \delta = \delta \ast f = f ,\;\;
\textrm{(4)}\;\; f \ast g = g \ast f$,

$\textrm{(5)}\;\; F( f \ast g )= (Ff)(Fg),\;\;
\textrm{(6)}\;\; \overline{F}( f \ast g )= (\overline{F} f)(\overline{F} g)$,

$\textrm{(7)}\;\; F(fg) = (Ff)\ast (Fg),\;\;
\textrm{(8)}\;\; \overline{F}(fg) = (\overline{F}f)\ast (\overline{F}g)$.

\medskip

The definition implies the following proposition :

\medskip

\textsc{Proposition} 1.5([N-O2]).\;\; If $l\in{\bf R}^{+}$, then
$F\delta^l=(H')^{(l-1)(\,^{\star}H)^2}$.

\medskip

If there exists $\alpha$, $\beta\in L^2({\bf R})$ so that
$a=\,^{\ast}\alpha|_L$, $b=\,^{\ast}\beta|_L$, that is,
$a(k)=\star(\,^{\ast}\alpha(k))$, $b(k)=\star(\,^{\ast}\beta(k))$, then
st(st($^{\star}\varepsilon\sum_{k\in
L}a(k)b(k)))=\int_{-\infty}^{\infty}a(x)b(x)dx$. 
Definition 1.3 is easier understanding than Definition 1.2 for a standard
meaning. For the reason, we consider mainly Definition 1.3 about several
examples.

\bigskip

{\bf 2. Examples of the infinitesimal Fourier transformation on the
space of functions}

\medskip

We calculate the infinitesimal Fourier transformations of
$\varphi_{\xi},\,\varphi_{im}\in R(L)$ :

\medskip

1. $\varphi_{\xi}(x)=\exp(-\xi\pi x^2),$ where $\xi\in{\bf C}$,
Re$(\xi)>0$, 

2. $\varphi_{im}(x)=\exp(-im\pi x^2)$,
where $m\in{\bf Z}$. 

\medskip

For $\varphi_{\xi}$, we obtain :

\medskip

{\bf Proposition
2.1.}

\medskip

$(F\varphi_{\xi})(p)=c_{\xi}(p)\varphi_{\xi}(\frac{p}{\xi})$,
where $c_{\xi}(p)=\sum_{x\in
L}\varepsilon\exp(-\xi\pi(x+\frac{i}{\xi}p)^2)$. 

\medskip
\noindent
If $p$ is finite, then
st$(c_{\xi}(p))=\frac{1}{\sqrt{\xi}}$. 

\medskip

{\bf Proof.}\;\;The infinitesimal Fourier
transformations of $\varphi_{\xi}$ is :

\medskip

$(F\varphi_{\xi})(p)=\sum_{x\in L}\varepsilon\exp(-2\pi ipx)\exp(-\xi\pi
x^2)$

$=\sum_{x\in
L}\varepsilon\exp(-\xi\pi(x+\frac{i}{\xi}p)^2-\pi\frac{1}{\xi}p^2)$

$=(\sum_{x\in
L}\varepsilon\exp(-\xi\pi(x+\frac{i}{\xi}p)^2))\exp(-\pi\frac{1}{\xi}p^2)=c_{\xi}(p)\varphi_{\xi}(\frac{p}{\xi})$,

\medskip
\noindent
where $c_{\xi}(p)=\sum_{x\in
L}\varepsilon\exp(-\xi\pi(x+\frac{i}{\xi}p)^2)$. If $p$ is finite, then
st$(c_{\xi}(p))\\=\int_{-\infty}^{\infty}\exp\left(-\xi\pi\left(t+\frac{i}{\xi}
\textrm{st}(p)\right)^2\right)dt=\frac{1}{\sqrt{\xi}}$. 

\medskip

Using Theorem 1.4(8), we obtain for $c_{\xi}$ :

\medskip

{\bf Proposition
2.2.}$\;\;\varphi_{\xi}(x')=\left(\overline{F}c_{\xi}(p)\ast
\left(c_{\frac{1}{\xi}}(-x)\varphi_{\xi}(x)\right)\right)(x')$.

\medskip

{\bf Proof.}\;\; We obtain :
$(F\varphi_{\xi})(p)=c_{\xi}(p)\varphi_{\xi}(\frac{p}{\xi})$, 
and put
$\overline{F}$ to the above :

\medskip

$(\overline{F}(F\varphi_{\xi}))(x)=(\overline{F}(c_{\xi}(p)\varphi_{\xi}
(\frac{p}{\xi})))(x)$

$=(\overline{F}c_{\xi}(p)\ast\overline
{F}\varphi_{\xi}(\frac{p}{\xi}))(x),$ that
is, $\varphi_{\xi}(x)=(\overline{F}c_{\xi}(p)\ast\overline
{F}\varphi_{\xi}(\frac{p}{\xi}))(x).$

\medskip
\noindent
Now $(\overline
{F}\varphi_{\xi}(\frac{p}{\xi}))(x)=\sum_{p\in
L}\varepsilon\exp(-2\pi ipx)\exp(-\xi(\frac{p}{\xi})^2\pi)$

$=\sum_{p\in
L}\varepsilon\exp(-\pi\frac{1}{\xi}(p^2-2\pi i\xi px))=\sum_{p\in
L}\varepsilon\exp(-\frac{\pi}{\xi}(p-i\xi x)^2+\frac{\pi}{\xi}(i\xi x)^2)$

$=\left(\sum_{p\in
L}\varepsilon\exp(-\frac{\pi}{\xi}(p-i\xi x)^2)\right)\exp(-\pi \xi
x^2)=\left(\sum_{p\in L}\varepsilon\exp(-\frac{\pi}{\xi}(p-i\xi
x)^2)\right)\varphi_{\xi}(x)$.

\medskip
\noindent
By the definition : $c_{\xi}(p)=\sum_{x\in
L}\varepsilon\exp(-\pi \xi(x+i\frac{1}{\xi}p)^2)$, the summation
\\$\sum_{p\in L}\varepsilon\exp(-\frac{\pi}{\xi}(p-i\xi
x)^2)$ is $c_{\frac{1}{\xi}}(-x)$.
Hence $\varphi_{\xi}(x')=\left(\overline{F}c_{\xi}(p)\ast
\left(c_{\frac{1}{\xi}}(-x)\varphi_{\xi}(x)\right)\right)(x')$.

\medskip

For the following proposition 2.3, we recall the Gauss sum(cf.[R]) :

\medskip

For $z\in {\bf N}$, Gauss sum
$\sum_{l=0}^{z-1}\exp(-i\frac{2\pi}{z}l^2)$ is equal to
$\sqrt{z}\frac{1+(-i)^z}{1-i}$.

\medskip

{\bf Proposition
2.3.}\;\;If $m|2H^2$ and $m|\frac {p}{\varepsilon}$, then
$(F\varphi_{im})(p)=c_{im}(p)\exp(i\pi\frac 1m p^2)$, where
$c_{im}(p)=\sqrt{\frac m2}\frac{1+i^{\frac{2H^2}{m}}}{1+i}$ for positive
$m$ and $c_{im}(p)=\sqrt{\frac
{-m}{2}}\frac{1+(-i)^{\frac{2H^2}{-m}}}{1-i}$ for negative
$m$.

\medskip

{\bf Proof.}\;\;$(F\varphi_{im})(p)=\sum_{x\in L}\varepsilon\exp(-im\pi
x^2)\exp(-2\pi ixp)$

$=\sum_{x\in L}\varepsilon\exp(-im\pi(x+\frac pm)^2)\exp(i\pi\frac 1m
p^2)$

$=c_{im}(p)\exp(i\pi\frac 1m p^2)$, where
$c_{im}(p)=\sum_{x\in L}\varepsilon\exp(-im\pi(x+\frac pm)^2)$.

\medskip
\noindent
Since $m|\frac {p}{\varepsilon}$, the element $\frac pm$ is in $L$. We
remark that $\exp(-i\pi mx^2)=\exp(-i\pi m(x+H)^2)$. For positive $m$,

\medskip

$c_{im}(p)=\sum_{x\in L}\varepsilon\exp(-im\pi x^2)=\sum_{0\le
n<H^2}\varepsilon\exp(-i2\pi\frac{m}{2H^2}n^2)$

\medskip

$=\frac m2 \sum_{0\le
n<\frac{2H^2}{m}}\varepsilon\exp(-i2\pi\frac{m}{2H^2}n^2)=\frac m2
\overline{\left(\varepsilon\sqrt{\frac{2H^2}{m}}
\frac{1+(-i)^{\frac{2H^2}{m}}}{1-i}\right)}$, by the above Gauss sum.
Hence $c_{im}(p)=\sqrt{\frac m2}\frac{1+i^{\frac{2H^2}{m}}}{1+i}$. For
negative $m$, the proof is as same as the above.

\bigskip 

{\bf 3. Examples of the infinitesimal Fourier transformation for the
space of functionals}

\medskip

We define an equivalence relation $\sim_{\,^{\star}HH'}$ in ${\bf
L}'$ by $x\sim_{\,^{\star}HH'} y\Leftrightarrow x-y\in
\,^{\star}HH'\,^{\star}(\,^{\ast}{\bf Z})$. We identify ${\bf
L}'/\sim_{\,^{\star}HH'}$ with
$L'$. Let

\medskip

$X_{H,\,^{\star}HH'}:=\{a'\,|\,a'\textrm{ is an internal function with a double meaning, from }\star({\bf L})/\sim_{\star(H)}\textrm{ to }\; {\bf
L}'/\sim_{\,^{\star}HH'}\}$, 

\medskip
\noindent
and let ${\bf e}$ be a mapping from $X$ to
$X_{H,\,^{\star}HH'}$, defined by $({\bf e}(a))([k])=[a(\hat k)]$, where
$[\;\;\;]$ in left hand side represents the equivalence class for the
equivalence relation $\sim_{\star(H)}$ in $\star({\bf L})$, $\hat k$ is a
representative in
$\star(L)$ satisfying $k \sim_{\star(H)} \hat k$, and $[\;\;\;]$ in right hand
side represents the equivalence class for the equivalence relation
$\sim_{\,^{\star}HH'}$ in ${\bf L}'$. Furthermore let ${\bf
e}^{\sharp}(f)(a')$ be defined by $f({\bf e}^{-1}(a'))$.

\medskip

{\bf 3-1. The infinitesimal Fourier transformation of $g_{\xi}(a)=\exp
\left(-\pi\,^{\star}\varepsilon\xi\sum_{k\in L}a^2(k)\right)$ with $\xi\in{\bf C}$,
Re$(\xi)>0$}

\medskip

We calculate the infinitesimal Fourier transformation
of

\medskip

$g_{\xi}(a)=\exp \left(-\pi\,^{\star}\varepsilon\xi\sum_{k\in 
L}a^2(k)\right)$, where $\xi\in {\bf C}$, Re$(\xi)>0$,

\medskip
\noindent
in the space $A$ of
functionals, for
Definition 1.3. We identify
$\,^{\star}(\,^{\ast}\xi)\in {\bf C}$ with
$\xi\in {\bf C}$. 

\medskip

{\bf Theorem 3.1.}$\;\;(F({\bf
e}^{\sharp}(g_{\xi})))(b)=C_{\xi}(b)g_{\xi}(\frac{b}{\xi}),$ where $b\in X$
and

$C_{\xi}(b)=\sum_{a\in X}\varepsilon_0\exp
\left(-\pi\,^{\star}\varepsilon\xi\sum_{k\in L}(a(k)+i\frac{1}{\xi}
b(k))^2\right).$

\medskip

{\bf Proof.}\;\; We do the infinitesimal Fourier transformation of
${\bf e}^{\sharp}(g_{\xi})(a)$.

\medskip
$(F({\bf e}^{\sharp}(g_{\xi})))(b)=F\left(\exp
\left(-\pi\,^{\star}\varepsilon\xi\sum_{k\in L}a^2(k)\right)\right)(b)$

$=\sum_{a\in X}\varepsilon_0\exp \left(-2i\pi
\,^{\star}\varepsilon\sum_{k\in L}a(k)b(k)\right)\exp
\left(-\pi\,^{\star}\varepsilon\xi\sum_{k\in L}a^2(k)\right)$

$=\sum_{a\in X}\varepsilon_0\exp
\left(-\pi\,^{\star}\varepsilon\xi\sum_{k\in L}(a^2(k)+2 i
\frac{1}{\xi}a(k)b(k))\right)$

$=\sum_{a\in
X}\varepsilon_0\exp\left(-\pi\,^{\star}\varepsilon\xi\sum_{k\in L}((a(k)+i
\frac{1}{\xi}b(k))^2+ \frac{1}{\xi}b^2(k))\right)$

$=\left(\sum_{a\in
X}\varepsilon_0\exp \left(-\pi\,^{\star}\varepsilon\xi\sum_{k\in
L}(a(k)+i \frac{1}{\xi}b(k))^2\right)\right)\exp \left(-\pi
\,^{\star}\varepsilon \frac{1}{\xi}\sum_{k\in L}b^2(k)\right)$

$=C_{\xi}(b)g_{\xi}(\frac{b}{\xi}).$

\medskip

Let $\star\circ\ast : {\bf R}\to \,^{\star}(\,^{\ast}{\bf R})$ be the natural
elementary embedding and let ${\bf st}(c)$ for $c\in\,^{\star}(\,^{\ast}{\bf
R})$ be the standard part of $c$ with respect to the natural elementary
embedding $\star\circ\ast$. 

\medskip

{\bf Theorem 3.2.}\;\;If the image of $b\;(\in X)$ is bounded by a finite
value of
$\,^{\ast}{\bf R}$, that is, $\exists b_0\in \,^{\ast}{\bf R}
\textrm{ s.t.
}k\in L\Rightarrow |b(k)|\le
\star(b_0),$
then

\medskip

$\textrm{st}(C_{\xi}(b))=\left(\ast\left(
\frac{1}{\sqrt{\xi}}\right)\right)^{H^2}\;(\in
\,^{\ast}{\bf R})$ and $\textrm{\bf
st}\left(\frac{C_{\xi}(b)}{\star\left(\left(\ast\left(\frac{1}{\sqrt{
\xi}}\right)\right)^{H^2}\right)}\right)=1.$

\medskip

{\bf Proof.}  We show
that
$\textrm{st}\left(\frac{C_{\xi}(b)}{\star\left(\left(\ast\left(\int_{-\infty}^{\infty}\exp
(-\pi \xi x^2)dx\right)\right)^{H^2}\right)}\right)=1.$

We consider the term $C_{\xi}(b)=\sum_{a\in X}\varepsilon_0\exp
\left(-\pi\,^{\star}\varepsilon\xi\sum_{k\in L}(a(k)+ i
\frac{1}{\xi}b(k))^2\right).$

\medskip
\noindent
We write
$(a(k))_{\lambda\mu}=\varepsilon_{\lambda\mu}'z_{\lambda\mu}^a
(k_{\mu})\;\;(z_{\lambda\mu}^a (k_{\mu})\in {\bf
Z}),\;(b(k))_{\lambda\mu}=\varepsilon_{\lambda\mu}'z_{\lambda\mu}^b
(k_{\mu})\;\;(z_{\lambda\mu}^b (k_{\mu})\in {\bf Z})$, where $\lambda\in\Lambda_1$,
$\mu\in\Lambda_2$. From now on, we
denote $z_{\lambda\mu}^a (k_{\mu})$, $z_{\lambda\mu}^b  (k_{\mu})$ by
$z_{\lambda\mu}^a$, $z_{\lambda\mu}^b$ for simplicity. Then the
$\lambda\mu$-component of $\sum_{a\in X}\varepsilon_0\exp 
\left(-\pi\,^{\star}\varepsilon\xi\sum_{k\in L}(a(k)+ i
\frac{1}{\xi}b(k))^2\right)$ is equal to

\medskip

$\sum_{a_{\lambda\mu}\in X_{\lambda\mu}}(\varepsilon_0)_{\lambda\mu} \exp
\left(-\pi\varepsilon_{\mu}\xi\sum_{k_{\mu}\in
L_{\mu}}((a(k))_{\lambda\mu}+ i
\frac{1}{\xi}(b(k))_{\lambda\mu})^2\right)$

$=\sum_{a_{\lambda\mu}\in
X_{\lambda\mu}}\prod_{k_{\mu}\in
L_{\mu}}\sqrt{\varepsilon_{\mu}}\varepsilon_{\lambda\mu}'\exp \left(-\pi
\varepsilon_{\mu}\xi((a(k))_{\lambda\mu}+ i
\frac{1}{\xi}(b(k))_{\lambda\mu})^2\right)$

$=\prod_{k_{\mu}\in
L_{\mu}}\left(\sum_{(a(k))_{\lambda\mu}\in
L_{\lambda\mu}'}\sqrt{\varepsilon_{\mu}}\varepsilon'{}_{\lambda\mu} \exp
\left(-\pi\varepsilon_{\mu}\xi ((a(k))_{\lambda\mu}+ i
\frac{1}{\xi}(b(k))_{\lambda\mu})^2\right)\right)$

$=\prod_{k_{\mu}\in
L_{\mu}}\left(\sum_{\varepsilon_{\lambda\mu}'z_{\lambda\mu}^a \in
L_{\lambda\mu}'}\sqrt{\varepsilon_{\mu}}\varepsilon'{}_{\lambda\mu} \exp
\left(-\pi\varepsilon_{\mu} \xi(\varepsilon_{\lambda\mu}'z_{\lambda\mu}^a+
i
\frac{1}{\xi}\varepsilon_{\lambda\mu}'z_{\lambda\mu}^b)^2\right)
\right)$

$=\prod_{k_{\mu}\in
L_{\mu}}\left(\sum_{\varepsilon_{\lambda\mu}'z_{\lambda\mu}^a \in
L_{\lambda\mu}'}\sqrt{\varepsilon_{\mu}}\varepsilon'{}_{\lambda\mu} \exp
\left(-\pi
\xi(\sqrt{\varepsilon_{\mu}}\varepsilon_{\lambda\mu}'z_{\lambda\mu}^a+ i
\frac{1}{\xi}\sqrt{\varepsilon_{\mu}}\varepsilon_{\lambda\mu}'
z_{\lambda\mu}^b)^2\right)\right).$

\medskip
\noindent
We assume that the image of $b\;(\in X)$ is bounded by a finite value of
$\,^{\ast}{\bf R}$, that is, $\exists b_0\in \,^{\ast}{\bf R}\textrm{ s.t.
}k\in L\Rightarrow |b(k)|\le
\star(b_0).$ The $\lambda\mu$-component of 

\noindent
$\frac{\sum_{a\in
X}\varepsilon_0\exp
\left(-\pi\,^{\star}\varepsilon\xi\sum_{k\in L}(a(k)+ i
\frac{1}{\xi}b(k))^2\right)}{\star\left(\left(\ast\left(
\int_{-\infty}^{\infty}\exp
(-\pi \xi x^2)dx\right)\right)^{H^2}\right)}$ is equal to

\noindent
$\prod_{k_{\mu}\in
L_{\mu}}\frac{\sum_{\varepsilon_{\lambda\mu}'z_{\lambda\mu}^a \in
L_{\lambda\mu}'}\sqrt{\varepsilon_{\mu}}\varepsilon'{}_{\lambda\mu} \exp
\left(-\pi
\xi(\sqrt{\varepsilon_{\mu}}\varepsilon_{\lambda\mu}'z_{\lambda\mu}^a+ i
\frac{1}{\xi}\sqrt{\varepsilon_{\mu}}\varepsilon_{\lambda\mu}'z_{\lambda\mu}^b)^2\right)}{\int_{-\infty}^{\infty}\exp(-\pi
\xi x^2)dx}.$ We write $\sqrt{\varepsilon_{\mu}}\varepsilon_{\lambda\mu}'z_{\lambda\mu}^a$,
$\sqrt{\varepsilon_{\mu}}\varepsilon_{\lambda\mu}'z_{\lambda\mu}^b$ as
$a_{\lambda\mu}$, $b_{\lambda\mu}$ for simplicity, and

\medskip

$B_{\lambda\mu}(k_{\mu}):=\sum_{\varepsilon_{\lambda\mu}'z_{\lambda\mu}^a
\in L_{\lambda\mu}'}\sqrt{\varepsilon_{\mu}}\varepsilon'{}_{\lambda\mu}
\exp \left(-\pi \xi(
\sqrt{\varepsilon_{\mu}}\varepsilon_{\lambda\mu}'z_{\lambda\mu}^a +
i\frac{1}{\xi}\sqrt{\varepsilon_{\mu}}\varepsilon_{\lambda\mu}'
z_{\lambda\mu}^b)^2\right)$

$-\int_{-\infty}^{\infty}\exp(-\pi
(x+i\frac{1}{\xi}\sqrt{\varepsilon_{\mu}}\varepsilon_{\lambda\mu}'
z_{\lambda\mu}^b)^2)dx.$

\medskip
\noindent
It is equal to

\medskip

$-\bigl(\int_{-\infty}^{-\sqrt{\varepsilon_{\mu}}\frac{H_{\lambda\mu}'}{2}}\exp(-\pi
\xi(x+ib_{\lambda\mu})^2)dx+\int_{\sqrt{\varepsilon_{\mu}}\frac{H_{\lambda\mu}'}{2}}^{\infty}\exp(-\pi
\xi(x+i\frac{1}{\xi} b_{\lambda\mu})^2)dx\bigr)$

$+\sum_{\varepsilon_{\lambda\mu}'z_{\lambda\mu}^a \in
L_{\lambda\mu}'}\sqrt{\varepsilon_{\mu}}\varepsilon'{}_{\lambda\mu} \exp
\left(-\pi \xi(
\sqrt{\varepsilon_{\mu}}\varepsilon_{\lambda\mu}'z_{\lambda\mu}^a+i
\frac{1}{\xi}\sqrt{\varepsilon_{\mu}}\varepsilon_{\lambda\mu}'
z_{\lambda\mu}^b)^2\right)$

$-\int_{-\sqrt{\varepsilon_{\mu}}\frac{H_{\lambda\mu}'}{2}}^{\sqrt{\varepsilon_{\mu}}\frac{H_{\lambda\mu}'}{2}}\exp(-\pi
\xi(x+i\frac{1}{\xi} b_{\lambda\mu})^2)dx\;\;\;\;\;\cdots (\ast_1).$

\medskip
\noindent
Then the above is equal to

\medskip

$\prod_{k_{\mu}\in
L_{\mu}}\frac{\sum_{\varepsilon_{\lambda\mu}'z_{\lambda\mu}^a \in
L_{\lambda\mu}'}\sqrt{\varepsilon_{\mu}}\varepsilon'{}_{\lambda\mu} \exp
\left(-\pi \xi(
\sqrt{\varepsilon_{\mu}}\varepsilon_{\lambda\mu}'z_{\lambda\mu}^a+ i 
\frac{1}{\xi}\sqrt{\varepsilon_{\mu}}\varepsilon_{\lambda\mu}'z_{\lambda\mu}^b)^2\right)}{\int_{-\infty}^{\infty}\exp(-\pi
\xi x^2)dx}$

$=\prod_{k_{\mu}\in
L_{\mu}}\left(1+\frac{B_{\lambda\mu}(k_{\mu})}{\int_{-\infty}^{\infty}\exp(-\pi
\xi x^2)dx}\right)
\;\;\;\;\;\cdots (\ast_2).$

\medskip

We show that $[(B_{\lambda\mu}(k_{\mu}))]$ is infinitesimal in 
$\,^{\star}(\,^{\ast}{\bf C})$ with respect to $\,^{\ast}{\bf C}$. It implies that
$\left[\left(\frac{1}{B_{\lambda\mu}(k_{\mu})}\right)\right]$ is infinite in
$\,^{\star}(\,^{\ast}{\bf C})$. Since
$b_{\lambda\mu}$ is finite and
$\left[\left(\sqrt{\varepsilon_{\mu}}\frac{H_{\lambda\mu}'}{2}
\right)\right]$ is infinitesimal in $\,^{\star}(\,^{\ast}{\bf
R})$ with respect to $\,^{\ast}{\bf R}$, the first and second terms 
of ($\ast_1$), that is,
$\left[\left(\int_{-\infty}^{-\sqrt{\varepsilon_{\mu}}\frac{H_{\lambda\mu}'}{2}}\exp(-\pi
\xi(x+i\frac{1}{\xi} b_{\lambda\mu})^2)dx\right)\right]$
and
$\left[\left(\int_{\sqrt{\varepsilon_{\mu}}\frac{H_{\lambda\mu}'}{2}}^{\infty}\exp(-\pi
\xi(x+i\frac{1}{\xi} b_{\lambda\mu})^2)dx\right)\right]$ is infinitesimal
in $\,^{\star}(\,^{\ast}{\bf C})$ with respect to $\,^{\ast}{\bf C}$. In order to
show that
$[(B_{\lambda\mu}(k_{\mu}))]$  is infinitesimal in $\,^{\star}(\,^{\ast}{\bf
C})$, we consider the third and forth terms in $(\ast_1)$, and we prove that it
is represents an infinitesimal number. 

First we calculate 

\medskip
$\exp(-\pi
\xi(x+i\frac{1}{\xi}\sqrt{\varepsilon_{\mu}}\varepsilon_{\lambda\mu}'z_{\lambda\mu}^b)^2)-\exp(-\pi
\xi(\sqrt{\varepsilon_{\mu}}\varepsilon_{\lambda\mu}'
z_{\lambda\mu}^a+i\frac{1}{\xi}\sqrt{\varepsilon_{\mu}}
\varepsilon_{\lambda\mu}'z_{\lambda\mu}^b)^2).$

\medskip
\noindent
Since

\medskip
$\exp(-\pi \xi(x+i\frac{1}{\xi}b_{\lambda\mu})^2)$

$=\exp(-\pi (\alpha x^2-\frac{\alpha
b_{\lambda\mu}^2}{\alpha^2+\beta^2}))\exp(-i\pi (\beta
x^2+2b_{\lambda\mu}x+\frac{\beta b_{\lambda\mu}^2}{\alpha^2+\beta^2}))$

$=\exp(-\pi (\alpha x^2-\frac{\alpha
b_{\lambda\mu}^2}{\alpha^2+\beta^2}))\cos(\pi (\beta
x^2+2b_{\lambda\mu}x+\frac{\beta b_{\lambda\mu}^2}{\alpha^2+\beta^2}))$

$-i\exp(-\pi (\alpha x^2-\frac{\alpha
b_{\lambda\mu}^2}{\alpha^2+\beta^2}))\sin(\pi (\beta
x^2+2b_{\lambda\mu}x+\frac{\beta b_{\lambda\mu}^2}{\alpha^2+\beta^2})),$

\medskip
\noindent
the above is

\medskip
$\exp(-\pi\xi (x+i\frac{1}{\xi}b_{\lambda\mu})^2)-\exp(-\pi\xi
(a_{\lambda\mu}+i\frac{1}{\xi}b_{\lambda\mu})^2)$

$=\bigl\{\exp(-\pi
(\alpha x^2-\frac{\alpha b_{\lambda\mu}^2}{\alpha^2+\beta^2}))\cos(\pi
(\beta x^2+2b_{\lambda\mu}x+\frac{\beta
b_{\lambda\mu}^2}{\alpha^2+\beta^2}))$

$-\exp(-\pi (\alpha
a_{\lambda\mu}^2-\frac{\alpha
b_{\lambda\mu}^2}{\alpha^2+\beta^2}))\cos(\pi (\beta
a_{\lambda\mu}^2+2b_{\lambda\mu}a_{\lambda\mu}+\frac{\beta
b_{\lambda\mu}^2}{\alpha^2+\beta^2}))\bigr\}$

$-i\bigl\{\exp(-\pi
(\alpha x^2-\frac{\alpha b_{\lambda\mu}^2}{\alpha^2+\beta^2}))\sin(\pi
(\beta x^2+2b_{\lambda\mu}x+\frac{\beta
b_{\lambda\mu}^2}{\alpha^2+\beta^2}))$

$-\exp(-\pi (\alpha
a_{\lambda\mu}^2-\frac{\alpha
b_{\lambda\mu}^2}{\alpha^2+\beta^2}))\sin(\pi (\beta
a_{\lambda\mu}^2+2b_{\lambda\mu}a_{\lambda\mu}+\frac{\beta
b_{\lambda\mu}^2}{\alpha^2+\beta^2}))\bigr\}.\;\;\;\;\;\cdots(\ast_3)$

\medskip
\noindent
We consider the first term of $(\ast_3)$. Then

\medskip
$\exp(-\pi (\alpha x^2-\frac{\alpha
b_{\lambda\mu}^2}{\alpha^2+\beta^2}))\cos(\pi (\beta
x^2+2b_{\lambda\mu}x+\frac{\beta b_{\lambda\mu}^2}{\alpha^2+\beta^2}))$

\medskip
\noindent
is equal to

\medskip

$\exp(\pi \frac{\alpha b_{\lambda\mu}^2}{\alpha^2+\beta^2})\exp(-\pi
\alpha x^2)\cos(\pi (\beta x^2+2b_{\lambda\mu}x+\frac{\beta
b_{\lambda\mu}^2}{\alpha^2+\beta^2})).$

\medskip
\noindent
We put

\medskip
$f(x)=\exp(-\pi \alpha x^2)\cos(\pi (\beta
x^2+2b_{\lambda\mu}x+\frac{\beta b_{\lambda\mu}^2}{\alpha^2+\beta^2})).$

\medskip
\noindent
We assume that $0\le b_{\lambda\mu}$.

\medskip

$f'(x)=-2\pi \alpha x\exp(-\pi \alpha x^2)\cos(\pi (\beta
x^2+2b_{\lambda\mu}x+\frac{\beta b_{\lambda\mu}^2}{\alpha^2+\beta^2}))$

$-\exp(-\pi \alpha x^2)(2\pi\beta x+2\pi b_{\lambda\mu})\sin(\pi (\beta
x^2+2b_{\lambda\mu}x+\frac{\beta b_{\lambda\mu}^2}{\alpha^2+\beta^2}))$

$=-2\pi\sqrt{(\alpha x)^2+(\beta x+ b_{\lambda\mu})^2}\exp(-\pi \alpha
x^2)\cos(\pi (\beta x^2+2b_{\lambda\mu}x+\frac{\beta
b_{\lambda\mu}^2}{\alpha^2+\beta^2})+\alpha_x),$

\medskip
\noindent
where

\medskip

$\cos\alpha_x=\frac{\alpha x}{\sqrt{(\alpha x)^2+(\beta x+
b_{\lambda\mu})^2}},\;\; -\sin \alpha_x=\frac{\beta
x+b_{\lambda\mu}}{\sqrt{(\alpha x)^2+(\beta x+ b_{\lambda\mu})^2}}.$

\medskip
\noindent
Since $0<\alpha$, if $0\le \beta$ and $0\le x$, then $0\le
\cos\alpha_x\textrm{ and }\sin \alpha_x\le 0.$
Hence $-\frac{\pi}{2}\le\alpha_x<0.$ There is a unique maximum of $|f(x)|$
in

\medskip

$\left\{x\in{\bf R}\,\left|\,\frac{\pi}{2}(2m-1)\le \pi (\beta
x^2+2b_{\lambda\mu}x+\frac{\beta
b_{\lambda\mu}^2}{\alpha^2+\beta^2})<\frac{\pi}{2}(2m+1)\right\}\right.$

\medskip
\noindent
for each $m\in{\bf Z}\; (-1\le m)$, that is, $x$ satisfies

\medskip

$f'(x)=0,\; \frac{\pi}{2}(2m-1)\le \pi (\beta
x^2+2b_{\lambda\mu}x+\frac{\beta
b_{\lambda\mu}^2}{\alpha^2+\beta^2})<\frac{\pi}{2}(2m+1)$

$\Longleftrightarrow \pi (\beta x^2+2b_{\lambda\mu}x+\frac{\beta
b_{\lambda\mu}^2}{\alpha^2+\beta^2})+\alpha_x=\frac{\pi}{2}(2m-1).\;\;\;\;\;\cdots
(\ast_4)$

\medskip
\noindent
We write the value of $x$ having the maximum of $|f(x)|$ in the interval as
$(A_{2m})_{\lambda\mu}$. On the other hand, we denote the value $\alpha_x$ at $x=(A_{2m})_{\lambda\mu}$  by
$\alpha_{A_{2m}}$. Then

$(A_{2m})_{\lambda\mu}=\frac{-b_{\lambda\mu}+\sqrt{\frac{\alpha^2b_{
\lambda\mu}^2}{\alpha^2+\beta^2}+\frac{\beta}{2}(2m-1)-\frac{
\alpha_{A_{2m}}}{\pi}}}{\beta}.$

\medskip
\noindent
The maximum of $f(x)$ is $f((A_{2m})_{\lambda\mu})=\exp\left(-\pi \alpha
(A_{2m})_{\lambda\mu}^2\right)\cos \left(m\pi -\frac{\pi}{2}
-\alpha_{A_{2m}} \right).$

\noindent
Since $\lim_{m\to\infty}\frac{(A_{2m})_{\lambda\mu}}
{\sqrt{\frac{2m-1}{2\beta}}}=1$, there exists $m$ such
that $\sqrt{\frac{2m-1}{4\beta}}<(A_{2m})_{\lambda\mu}.$ Hence there exists
$m$ such that 

\medskip
$|f((A_{2m})_{\lambda\mu})|\le \exp\left(-\pi
(A_{2m})_{\lambda\mu}^2\right)\le \exp\left(-\pi \frac{2m-1}{4\beta}\right).$

\medskip
\noindent
We denote the value of $x$ at $f(x)=0$, that is,

\medskip
$\pi (\beta x^2+2b_{\lambda\mu}x+\frac{\beta
b_{\lambda\mu}^2}{\alpha^2+\beta^2})=\frac {\pi}{2}(2m+1)$

\medskip
\noindent
by $(A_{2m+1})_{\lambda\mu}$.
Then
$(A_{2m+1})_{\lambda\mu}=\frac{-b_{\lambda\mu}+\sqrt{\frac{
\alpha^2b_{\lambda\mu}^2}{\alpha^2+\beta^2}+\frac{\beta}{2}(2m-1)}}
{\beta}.$ We consider

\medskip

$\Bigl|\sum_{\varepsilon_{\lambda\mu}'z_{\lambda\mu}^a \in
L_{\lambda\mu}'}\sqrt{\varepsilon_{\mu}}\varepsilon'{}_{\lambda\mu}
\exp(-\pi
(\sqrt{\varepsilon_{\mu}}\varepsilon_{\lambda\mu}'z_{\lambda\mu}^a)^2-(\varepsilon_{\lambda\mu}'z_{\lambda\mu}^b)^2))
\cos(2\pi
\varepsilon_{\lambda\mu}'z_{\lambda\mu}^b\sqrt{\varepsilon_{\mu}}
\varepsilon_{\lambda\mu}'z_{\lambda\mu}^a)
-\int_{-\sqrt{\varepsilon_{\mu}}\frac{H_{\lambda\mu}'}{2}}^{\sqrt{\varepsilon_{\mu}}\frac{H_{\lambda\mu}'}{2}}\exp(-\pi
(x^2-b_{\lambda\mu}^2))\cos(2\pi b_{\lambda\mu}x)dx\Bigr|.$

\medskip
\noindent
It is equal to

\medskip

$\exp(\pi
b_{\lambda\mu}^2)\Bigl|\sum_{\varepsilon_{\lambda\mu}'z_{\lambda\mu}^a
\in L_{\lambda\mu}'}\sqrt{\varepsilon_{\mu}}\varepsilon'{}_{\lambda\mu}
\exp(-\pi
(\sqrt{\varepsilon_{\mu}}\varepsilon_{\lambda\mu}'z_{\lambda\mu}^a)^2)
\cos(2\pi \varepsilon_{\lambda\mu}'z_{\lambda\mu}^b\sqrt{
\varepsilon_{\mu}}\varepsilon_{\lambda\mu}'z_{\lambda\mu}^a)-\int_{-\sqrt{\varepsilon_{\mu}}\frac{H_{\lambda\mu}'}{2}}^{\sqrt{\varepsilon_{\mu}}\frac{H_{\lambda\mu}'}{2}}\exp(-\pi
x^2)\cos(2\pi b_{\lambda\mu}x)dx\Bigr|$

$=\exp(\pi
b_{\lambda\mu}^2)\Bigl|\sum_{\varepsilon_{\lambda\mu}'z_{\lambda\mu}^a
\in L_{\lambda\mu}'}\sqrt{\varepsilon_{\mu}}\varepsilon'{}_{\lambda\mu}
f(\sqrt{\varepsilon_{\mu}}\varepsilon_{\lambda\mu}'z_{\lambda\mu}^a)^2)
-\int_{-\sqrt{\varepsilon_{\mu}}\frac{H_{\lambda\mu}'}{2}}^{\sqrt
{\varepsilon_{\mu}}\frac{H_{\lambda\mu}'}{2}}f(x)dx\Bigr|.$

\medskip
\noindent
We show that the following term is
infinitesimal in $\,^{\star}(\,^{\ast}{\bf R})$ with respect to $\,^{\ast}{\bf R}$ :

\medskip

$\Bigl|\sum_{\varepsilon_{\lambda\mu}'z_{\lambda\mu}^a \in
L_{\lambda\mu}'}\sqrt{\varepsilon_{\mu}}\varepsilon'{}_{\lambda\mu}
f(\sqrt{\varepsilon_{\mu}}\varepsilon_{\lambda\mu}'z_{
\lambda\mu}^a)^2)-\int_{-\sqrt{\varepsilon_{\mu}}\frac{H_{\lambda\mu}'}
{2}}^{\sqrt{\varepsilon_{\mu}}\frac{H_{\lambda\mu}'}{2}}f(x)dx\Bigr|.$

\medskip

Now

\medskip

$(A_{2(m+1)+1})_{\lambda\mu}-(A_{2m+1})_{\lambda\mu}$

$=
\frac{-b_{\lambda\mu}+\sqrt{\frac{\alpha^2b_{\lambda\mu}^2}{\alpha^2+\beta^2}+\frac{\beta}{2}(2(m+1)-1)}}{\beta}
-\frac{-b_{\lambda\mu}+\sqrt{\frac{\alpha^2b_{\lambda\mu}^2}{\alpha^2+
\beta^2}+\frac{\beta}{2}(2m-1)}}{\beta}$

$=\frac{\sqrt{\frac{\alpha^2b_{\lambda\mu}^2}{\alpha^2+\beta^2}+\frac{\beta}
{2}(2(m+1)-1)}-\sqrt{\frac{\alpha^2b_{\lambda\mu}^2}{\alpha^2+\beta^2}
+\frac{\beta}{2}(2m-1)}}{\beta}$

$=\frac{1}{2\left(\sqrt{\frac{\alpha^2b_{\lambda\mu}^2}{\alpha^2+
\beta^2}+\frac{\beta}{2}(2(m+1)-1)}+\sqrt{\frac{\alpha^2b_{\lambda\mu}
^2}{\alpha^2+\beta^2}+\frac{\beta}{2}(2m-1)}\right)}\;\;\;\;\;\cdots
(\ast_5).$

\medskip
\noindent
Since 

\medskip

$-\sqrt{\varepsilon_{\mu}}\frac{H_{\lambda\mu}'}{2}\le x\le
\sqrt{\varepsilon_{\mu}}\frac{H_{\lambda\mu}'}{2}
\;\;\;\;\;\cdots(\ast_6)$

\medskip
\noindent
and the image of $b\;(\in X)$ is bounded by a finite value of $\,^{\ast}{\bf R}$, that
is, $\exists b_0\in \,^{\ast}{\bf R}\textrm{ s.t. }k\in L\Rightarrow
|b(k)|\le
\star(b_0),$ 
the above ($\ast_5$) is greater than the following value
: $\frac{1}{4\sqrt{|b_0|^2+\frac{\beta}{2}(2m+1)}}.$ The value $\pi(2m+1)$
satisfies ($\ast_6$), $\pi(2m+1)\le
\sqrt{\varepsilon_{\mu}}\frac{H_{\lambda\mu}'}{2},$ that is, $2m+1\le
\frac{\sqrt{\varepsilon_{\mu}}H_{\lambda\mu}'}{2\pi},$
and $\frac{1}{4\sqrt{|b_0|^2+\frac{2}{\beta}(2m+1)}}\geq
\frac{1}{4\sqrt{|b_0|^2+\frac{\beta\sqrt{\varepsilon_{\mu}}
H_{\lambda\mu}'}{4\pi}}}.$ 

\noindent
Furthermore $\frac{1}{4\sqrt{|b_0|^2+\frac{2}{\beta}(2m+1)}}\geq
\frac{1}{4\sqrt{|b_0|^2+\frac{\beta\sqrt{\varepsilon}H'}{4\pi}}}>
\sqrt{\varepsilon}\varepsilon'.$

\noindent
Hence

\medskip

$\left\{\lambda\,\left|\,\left\{\mu\,\left|\,|(A_{2(m+1)+1})_{\lambda\mu}-(A_{2m+1})_{\lambda\mu}|>\sqrt{\varepsilon_{\lambda}}\varepsilon_{\lambda\mu}'\right\}\in
F_2\right\}\in F_1\right.\right. .$

\medskip
\noindent
We consider the $\lambda\mu$-component satisfying the above. Now

\medskip

$\Bigl|\sum_{\varepsilon_{\lambda\mu}'z_{\lambda\mu}^a \in
L_{\lambda\mu}'}\sqrt{\varepsilon_{\mu}}\varepsilon'{}_{\lambda\mu}
f(\sqrt{\varepsilon_{\mu}}\varepsilon_{\lambda\mu}'z_{\lambda\mu}^a)^2)
-\int_{-\sqrt{\varepsilon_{\mu}}\frac{H_{\lambda\mu}'}{2}}^
{\sqrt{\varepsilon_{\mu}}\frac{H_{\lambda\mu}'}{2}}f(x)dx\Bigr|$

$\le \Bigl|\sum_{\varepsilon_{\lambda\mu}'z_{\lambda\mu}^a \in
L_{\lambda\mu}',\; z_{\lambda\mu}^a \geq 0}
\sqrt{\varepsilon_{\mu}}\varepsilon'{}_{\lambda\mu} 
f(\sqrt{\varepsilon_{\mu}}\varepsilon_{\lambda\mu}'z_{\lambda\mu}^a)
^2)-\int_{0}^{\sqrt{\varepsilon_{\mu}}\frac{H_{\lambda\mu}'}{2}}f(x)dx
\Bigr|$

$+\Bigl|\sum_{\varepsilon_{\lambda\mu}'z_{\lambda\mu}^a \in
L_{\lambda\mu}',\; z_{\lambda\mu}^a\le
0}\sqrt{\varepsilon_{\mu}}\varepsilon'{}_{\lambda\mu}
f(\sqrt{\varepsilon_{\mu}}\varepsilon_{\lambda\mu}'z_{\lambda\mu}^a)
^2)-\int_{-\sqrt{\varepsilon_{\mu}}\frac{H_{\lambda\mu}'}{2}}^{0}f(x)
dx\Bigr|.$

\medskip
\noindent
We consider the difference

\medskip

$\sum_{\varepsilon_{\lambda\mu}'z_{\lambda\mu}^a \in L_{\lambda\mu}',\; z_{\lambda\mu}^a \geq 0}
\sqrt{\varepsilon_{\mu}}\varepsilon'{}_{\lambda\mu} f(\sqrt{\varepsilon_{\mu}}\varepsilon_{\lambda\mu}'z_{\lambda\mu}^a)^2)-\int_{0}^{\sqrt{\varepsilon_{\mu}}\frac{H_{\lambda\mu}'}{2}}f(x)dx.$

\medskip
\noindent
We devide the interval $(-\infty, \infty)$ into a sum of intervals where the 
function $f$ is monoton increasing or monoton decreasing. The absolute value of the
difference is bounded to the sum of the absolute values of the difference whose integral
areas are restricted to these intervals. Each difference is bounded to the product of
$\{$(the maximum value of $f$ on the interval) $-$ (the minimum value of $f$ on the
interval)$\}$ and
$\sqrt{\varepsilon_{\mu}}\varepsilon'{}_{\lambda\mu}$. Hence

\medskip
$\Bigl|\sum_{\varepsilon_{\lambda\mu}'z_{\lambda\mu}^a \in
L_{\lambda\mu}',\; z_{\lambda\mu}^a \geq 0}
\sqrt{\varepsilon_{\mu}}\varepsilon'{}_{\lambda\mu} 
f(\sqrt{\varepsilon_{\mu}}\varepsilon_{\lambda\mu}'z_{\lambda\mu}
^a)^2)-\int_{0}^{\sqrt{\varepsilon_{\mu}}\frac{H_{\lambda\mu}'}{2}}
f(x)dx\Bigr|$

$\le\sqrt{\varepsilon_{\mu}}\varepsilon'{}_{\lambda\mu}
2\sum_{m=0}^{\infty}|f(A_{2m})|$

$\le\sqrt{\varepsilon_{\mu}}\varepsilon'{}_{\lambda\mu}
2\sum_{m=0}^{\infty}\exp\left(-\pi
\left(\frac{2m-1}{4\beta}\right)\right).$

\medskip
\noindent
Since the value $2\sum_{m=0}^{\infty}\exp\left(-\pi
\left(\frac{2m-1}{4\beta}\right)\right)$ is finite, the following 
value

\medskip
$\sqrt{\varepsilon_{\mu}}\varepsilon'{}_{\lambda\mu}
2\sum_{m=0}^{\infty}\exp\left(-\pi
\left(\frac{2m-1}{4\beta}\right)\right)$

\medskip
\noindent
is infinitesimal in $\,^{\star}(\,^{\ast}{\bf R})$ with respect to 
$\,^{\ast}{\bf R}$.

\medskip

The same argument implies in the case $x<0$

\medskip

$\Bigl|\sum_{\varepsilon_{\lambda\mu}'z_{\lambda\mu}^a \in
L_{\lambda\mu}',\; z_{\lambda\mu}^a\le
0}\sqrt{\varepsilon_{\mu}}\varepsilon'{}_{\lambda\mu}
f(\sqrt{\varepsilon_{\mu}}\varepsilon_{\lambda\mu}'z_{\lambda\mu}^a)
^2)-\int_{-\sqrt{\varepsilon_{\mu}}\frac{H_{\lambda\mu}'}{2}}^{0}f(x)
dx\Bigr|$

\medskip
\noindent
is infinitesimal in $\,^{\star}(\,^{\ast}{\bf R})$ with respect to 
$\,^{\ast}{\bf R}$ also.

Hence

\medskip

$\Bigl|\sum_{\varepsilon_{\lambda\mu}'z_{\lambda\mu}^a \in
L_{\lambda\mu}'}\sqrt{\varepsilon_{\mu}}\varepsilon'{}_{\lambda\mu}
f(\sqrt{\varepsilon_{\mu}}\varepsilon_{\lambda\mu}'z_{\lambda\mu}^a)
^2)-\int_{-\sqrt{\varepsilon_{\mu}}\frac{H_{\lambda\mu}'}{2}}^{\sqrt
{\varepsilon_{\mu}}\frac{H_{\lambda\mu}'}{2}}f(x)dx\Bigr|$

\medskip
\noindent
is infinitesimal in $\,^{\star}(\,^{\ast}{\bf R})$ with respect to $\,^{\ast}{\bf R}$.

If $b_{\lambda\mu}<0$, 
the argument is parallel, and also, for the term of sin in $(\ast_3)$, though
sin is not an even function, the same argument holds. Hence
$[(B_{\lambda\mu}(k_\mu))]$ is infinitesimal in $\,^{\star}(\,^{\ast}{\bf C})$
with respect to $\,^{\ast}{\bf C}$. Hence

\medskip

$\textrm{st}\left(\frac{\sum_{a\in X}\varepsilon_0\exp
\left(-\pi\,^{\star}\varepsilon\xi \sum_{k\in L}(a(k)+i\frac{1}{\xi}
b(k))^2\right)}{\star\left(\left(\ast\left(\int_{-\infty}^{\infty}\exp
(-\pi \xi x^2)dx\right)\right)^{H^2}\right)}\right)$

$=\textrm{st}\left(\prod_{k_{\mu}\in
L_{\mu}}\left(1+\frac{B_{\lambda\mu}(k_{\mu})}{\int_{-\infty}^{\infty}\exp(-\pi
\xi x^2)dx}\right)\right)$

$=\prod_{k_{\mu}\in
L_{\mu}}\textrm{st}\left(\left(1+\frac{B_{\lambda\mu}(k_{\mu})}{\int_{-\infty}^{\infty}\exp(-\pi
\xi x^2)dx}\right)\right)$

$=\prod_{k_{\mu}\in L_{\mu}}1$

$=1.$

\medskip
\noindent
Since $\int_{-\infty}^{\infty}\exp (-\pi \xi x^2)dx=\frac{1}
{\sqrt{\xi}},$ then 

\medskip

$\textrm{st}\left(\frac{\sum_{a\in X}
\varepsilon_0\exp
\left(-\pi\,^{\star}\varepsilon\xi \sum_{k\in L}(a(k)+i\frac{1}{\xi}
b(k))^2\right)}{\star\left(\left(\ast\left(\frac{1}{\sqrt{\xi}}\right)
\right)^{H^2}\right)}\right)=1,$ that
is, $\textrm{st}\left(\frac{C_{\xi}(b)}{\star\left(\left(\ast\left(
\frac{1}{\sqrt{\xi}}\right)\right)^{H^2}\right)}\right)=1.$

\medskip
\noindent
Furthermore

\medskip

$\textrm{\bf
st}\left(\frac{C_{\xi}(b)}{\star\left(\left(\ast\left(\frac{1}{\sqrt{
\xi}}\right)\right)^{H^2}\right)}\right)=\textrm{st}\left(\textrm{st}
\left(\frac{C_{\xi}(b)}{\star\left(\left(\ast\left(\frac{1}{\sqrt{\xi}}
\right)\right)^{H^2}\right)}\right)\right)
=1.$

\medskip

The argument is same about the infinitesimal Fourier
transformation of $g_{\xi}'(a)=\exp(-\pi\xi\sum_{k\in L}a^2(k))$, for
Definition 1.2, as the above.

\medskip

{\bf Theorem 3.3.}$\;\;(F({\bf
e}^{\sharp}(g_{\xi}')))(b)=B_{\xi}(b)g_{\xi}'(\frac{b}{\xi}),$ where $b\in
X$ and

$B_{\xi}(b)=\sum_{a\in X}\varepsilon_0\exp
\left(-\pi\xi\sum_{k\in L}(a(k)+i\frac{1}{\xi}
b(k))^2\right)$. Furthermore, if the image of $b\;(\in X)$ is bounded by a finite
value of
$\,^{\ast}{\bf R}$, that is, $\exists b_0\in \,^{\ast}{\bf R}
\textrm{ s.t.
}k\in L\Rightarrow |b(k)|\le
\star(b_0),$
then

\medskip

$\textrm{st}(B_{\xi}(b))=\left(\ast\left(
\frac{1}{\sqrt{\xi}}\right)\right)^{H^2}\;(\in
\,^{\ast}{\bf R})$ and $\textrm{\bf
st}\left(\frac{B_{\xi}(b)}{\star\left(\left(\ast\left(\frac{1}{\sqrt{
\xi}}\right)\right)^{H^2}\right)}\right)=1.$

\bigskip

{\bf 3-2. The infinitesimal Fourier transformation of $g_{im}=\exp(-i\pi
m\,^{\star}\varepsilon\sum_{k\in L}a^2(k))$ with $m\in {\bf Z}$}

\medskip

We calculate the infinitesimal Fourier transformation of

\medskip

$g_{im}(a)=\exp(-i\pi m\,^{\star}\varepsilon\sum_{k\in L}a^2(k)),$
where $m\in {\bf Z}$, 

\medskip
\noindent
for Definition 1.3.

\medskip

{\bf Proposition 3.4.}\;\;$(F({\bf e}^{\sharp}(g_{im})))(b)$
is written as $C_{im}(b)g_{\frac{1}{im}}(b)$. 

\medskip
\noindent
If $m|2\,^{\star}HH'{}^2$
and
$m|\frac{b(k)}{\varepsilon'}$ for an arbitrary $k$ in $L$,
then $(F({\bf e}^{\sharp}(g_{im})))(b)=C_{im}(b)g_{\frac{1}{im}}(b)$,
where $C_{im}(b)=\left(\sqrt{\frac
m2}\frac{1+i^{\frac{2\,^{\star}HH'{}^2}{m}}}{1+i}\right)^
{(\,^{\star}H)^2}$
for a positive $m$ and 

\noindent
$C_{im}(b)=\left(\sqrt{\frac
{-m}{2}}\frac{1+(-i)^{\frac{2\,^{\star}HH'{}^2}{-m}}}{1-i}\right)
^{(\,^{\star}H)^2}$
for a negative $m$.

\medskip

{\bf Proof.}$\;\;(F({\bf e}^{\sharp}(g_{im})))(b)=\sum_{a\in
X}\varepsilon_0\exp(-2\pi i\,^{\star}\varepsilon\sum_{k\in
L}a(k)b(k))\exp(-i\pi m
\,^{\star}\varepsilon\sum_{k\in L}a^2(k))$

$=\sum_{a\in X}\varepsilon_0\exp(-i\pi m
\,^{\star}\varepsilon\sum_{k\in L}(a(k)+\frac 1m b(k))^2)\exp(i\pi \frac
1m
\,^{\star}\varepsilon\sum_{k\in L}b^2(k))$

$=C_{im}(b)g_{\frac{1}{im}}(b)$,
where $C_{im}(b)=\sum_{a\in X}\varepsilon_0\exp(-i\pi m
\,^{\star}\varepsilon\sum_{k\in L}(a(k)+\frac 1m b(k))^2).$ 

\medskip
\noindent
When we denote
$a(k)$, $b(k)$ by $\varepsilon'n'$, $\varepsilon'l'$, 

\medskip

$\sum_{-\,^{\star}H\frac{H'{}^2}{2}\le a(k)<\,^{\star}H\frac{H'{}^2}{2}
}\exp(-i\pi m
\,^{\star}\varepsilon\sum_{k\in L}(a(k)+\frac 1m b(k))^2)$

$=\sum_{-\,^{\star}H\frac{H'{}^2}{2}\le a(k)<\,^{\star}H\frac{H'{}^2}{2}
}\exp(-i\pi m
\,^{\star}\varepsilon\sum_{k\in
L}(\varepsilon'n'+\varepsilon'\frac{n'}{m})^2)$.

\medskip
\noindent
Since $m|\frac{b(k)}{\varepsilon'}$, for
positive $m$, it is equal to

\medskip

$\sum_{-\,^{\star}H\frac{H'{}^2}{2}\le a(k)<\,^{\star}H\frac{H'{}^2}{2}
}\exp(-i\pi m
\,^{\star}\varepsilon\varepsilon'{}^2
n'{}^2)=\sum_{-\,^{\star}H\frac{H'{}^2}{2}\le
a(k)<\,^{\star}H\frac{H'{}^2}{2} }\exp(-2\pi
i\frac{m}{2\,^{\star}HH'{}^2}n'{}^2)$

$=\frac{m}{2}\sum_{0\le n'<\frac{2\,^{\star}HH'{}^2}{m}}\exp(-2\pi
i\frac{m}{2\,^{\star}HH'{}^2}n'{}^2)=\frac{m}{2}\sqrt{\frac{2\,^{\star}HH'{}^2}{m}}\frac{1+i^{\frac{2\,^{\star}
HH'{}^2}{m}}}{1+i}$, 

\medskip
\noindent
by Lemma 2.3. Hence
$C_{im}=\left(\sqrt{\frac{m}{2}}\frac{1+i^{\frac{2\,^{\star}
HH'{}^2}{m}}}{1+i}\right)^{(\,^{\star}H)^2}$ for a positive $m$. For
a negative $m$, the proof is as same as the above.

\medskip

The argument is same about the infinitesimal Fourier
transformation of 

\medskip

$g_{im}'(a)=\exp(-i\pi m\sum_{k\in
L}a^2(k)),$ where $m\in {\bf Z}$, 

\medskip
\noindent
for
Definition 1.2, as the above.

\medskip

{\bf Proposition 3.5.}$\;\;$If $m|2\,^{\star}HH'{}^2$ and
$m|\frac{b(k)}{\varepsilon'}$ for an arbitrary $k$ in $L$,
then $(F({\bf e}^{\sharp}(g_{im}')))(b)=B_{im}(b)g_{\frac{1}{im}}'(b)$,
where $B_{im}(b)=\left(\sqrt{\frac
m2}\frac{1+i^{\frac{2H'{}^2}{m}}}{1+i}\right)^
{(\,^{\star}H)^2}$
for a positive $m$ and 
$B_{im}(b)=\left(\sqrt{\frac
{-m}{2}}\frac{1+(-i)^{\frac{2H'{}^2}{-m}}}{1-i}\right)
^{(\,^{\star}H)^2}$
for a negative $m$.

\bigskip

{\bf 4.  Poisson summation formula for infinitesimal Fourier
transformation by Kinoshita}

\medskip

We extend Poisson summation formula of finite group to Kinoshita's infinitesimal Fourier transformation. 

\medskip

{\bf 4-1.  Formulation}

\medskip

{\bf Theorem 4.1.}  Let $S$ be an internal subgroup of $L$. Then we
obtain, for $\varphi\in R(L)$, 

\medskip

$|S^\bot|^{-\frac 12}\sum_{p\in
S^\bot}(F\varphi)(p)=|S|^{-\frac 12}\sum_{x\in S}\varphi(x),$

\medskip
\noindent
where $S^\bot :=\{p\in L\,|\,
\exp(2\pi ipx)=1\textrm{ for }\forall x\in S\}.$

\medskip

Since $L$ is an internal cyclic group, $S$ is also  an internal cyclic 
group. The genarater of $L$ is $\varepsilon$. The genarater of $S$ is
written as $\varepsilon s\; (s\in \,^{\ast}{\bf Z})$. Since the order of
$L$ is $H^2$, so $s$ is a factor of $H^2$.

We prepare the following lemma for the proof of Theorem 4.1.

\medskip

{\bf Lemma 4.2.}\;\;$S^\bot =<\varepsilon\frac{H^2}{s}>$.

\medskip

{\bf Proof of Lemma 4.2.}\;\;For $p\in S^\bot$, we write $p=\varepsilon
t$. Then we obtain the following :

\medskip

$\exp(2\pi ip\varepsilon s)=1\Longleftrightarrow \exp(2\pi i\varepsilon
t\varepsilon s)=1\Longleftrightarrow \exp(2\pi
t\frac{s}{H^2})=1\Longleftrightarrow t\frac{s}{H^2}\in\,^{\ast}{\bf Z}.$

\medskip
\noindent
Hence the generater of $S^\bot$ is $\varepsilon\frac{H^2}{s}$.

\medskip

{\bf Proof of Theorem 4.1.}\;\;By Lemma 4.2, $|S|=\frac{H^2}{s}$ and
$|S^\bot|=s$. If $x\notin S$, then
$\varepsilon\frac{H^2}{s}xs=\varepsilon H^2 x\in \,^{\ast}{\bf
Z},\;\left(\exp\left(2\pi i\varepsilon\frac{H^2}{s}x\right)\right)^s=1.$
For $x\in L$,

\begin{align*}
\sum_{p\in S^\bot}\exp(2\pi ipx)=&
\begin{cases}
\frac{\exp(2\pi (-\frac H2)x)(1-(\exp(2\pi i\varepsilon\frac{H^2}
{s}x)^s))}{1-\exp(2\pi i\varepsilon\frac{H^2}{s}x)}
& \text{$(x\notin S)$}\\
\sum_{p\in S^\bot}1  & \text{$(x\in S)$}
\end{cases}\\
=&
\begin{cases}
0 & \text{$(x\notin S)$}\\
s  & \text{$(x\in S)$}
\end{cases}.
\end{align*}

\noindent
Hence

\medskip

$\sum_{p\in S^\bot}(F\varphi)(p)=\sum_{p\in S^\bot}\varepsilon(\sum_{x\in
L}\varphi(x)\exp(2\pi i px))=\varepsilon\sum_{x\in L}\varphi(x)(\sum_{p\in
S^\bot}\exp(2\pi i px))$

$=\varepsilon\sum_{x\in S}\varphi(x) s
=\frac sH\sum_{x\in S}\varphi(x).
$

\medskip
\noindent
Thus

$\frac{1}{\sqrt{s}}\sum_{p\in S^\bot}(F\varphi)(p)=\frac{1}{\sqrt{s}}\frac
sH\sum_{x\in S}\varphi(x) =\sqrt{\frac {s}{H^2}}\sum_{x\in
S}\varphi(x)\;\;\;\;\;\cdots (1),$

$|S^\bot|^{-\frac 12}\sum_{p\in S^\bot}(F\varphi)(p)=\frac{1}{|S|^{\frac
12}}\sum_{x\in S}\varphi(x).$

\medskip

{\bf Proposition 4.3}\;\; Especially if $s$ is equal to $H$, then $(1)$
implies that

\medskip

$\sum_{p\in S^\bot}(F\varphi)(p)=\sum_{x\in S}\varphi(x)$.

\medskip
\noindent
The standard part of the above is 

\medskip

st$(\sum_{p\in S^\bot}
(F\varphi)(p))=$st$(\sum_{x\in S}\varphi(x))$.

\medskip
\noindent
If there exists a standard function $\varphi':{\bf R}\to {\bf C}$ so that
$\varphi=\,^{\ast}\varphi'|_L$, then the right hand side is equal to
$\sum_{-\infty <x<\infty}\varphi'(x)$, that is, $\sum_{-\infty
<x<\infty}$st$(\varphi)(x)$. Furthermore if $\varepsilon s$ is
infinitesimal and $\varphi'$ is integrable on ${\bf R}$, then 

\medskip

st$(\varepsilon s \sum_{x\in
S}\varphi(x))=\int_{-\infty}^{\infty}\varphi'(x)dx$.

\medskip
\noindent
Since $(1)$ implies
that 

\medskip

$\sum_{p\in S^\bot}(F\varphi)(p)=\varepsilon s \sum_{x\in
S}\varphi(x)$, 

\medskip
\noindent
we obtain st$(\sum_{p\in
S^\bot}(F\varphi)(p))=\int_{-\infty}^{\infty}\varphi'(x)dx$, that is,
$\int_{-\infty}^{\infty}$st$(\varphi)(x)dx$.

\medskip

We decompose $H$ to prime factors $H=p_1^{l_1}p_2^{l_2}\cdots p_m^{l_m}$,
where $p_1=2,\; p_1<p_2<\cdots <p_m$, each $p_i$ is a prime number,
$0<l_i$. Since $S$ is a subgroup of $L$, the number $s$ is a factor of
$H^2$. When we write $s$ as $p_1^{k_1}p_2^{k_2}\cdots p_m^{k_m}$, the
order of $S$ is equal to $p_1^{2l_1-k_1}p_2^{2l_2-k_2}\cdots
p_m^{2l_m-k_m}$ and the order of $S^\bot$ is $p_1^{k_1}p_2^{k_2}\cdots
p_m^{k_m}$. Hence $(1)$ is 

\medskip

$(p_1^{k_1}p_2^{k_2}\cdots p_m^{k_m})^{-\frac 12}\sum(_{p\in
S^\bot}(F\varphi)(p))=(p_1^{2l_1-k_1}p_2^{2l_2-k_2}\cdots
p_m^{2l_m-k_m})^{-\frac 12}\sum_{x\in
S}\varphi(x)$.

\bigskip

{\bf 4-2.  Examples}

\medskip

We apply Theorem 4.1 to the following two functions :

\medskip

1. $\varphi_i(x)=\exp(-i\pi x^2),$

2. $\varphi_{\xi}(x)=\exp(-\xi\pi x^2).$

\medskip
\noindent
The infinitesimal Fourier transformations of the functions are :

\medskip

1. $(F\varphi_i)(p)=\exp(-i\frac{\pi}{4})\overline{\varphi_i(p)}$,

2. $(F\varphi_{\xi})(p)=c_{\xi}(p){\varphi_{\xi}(\frac{p}{\xi})}$,

\medskip
\noindent
where st$(c_{\xi}(p))=\frac{1}{\sqrt{\xi}}$, if $p$ is finite.
Hence we obtain :

\medskip

1. $|S^{\bot}|^{-\frac 12}\exp(-i\frac{\pi}{4})\sum_{p\in
S^{\bot}}\overline{\varphi_i(p)}=|S|^{-\frac 12}\sum_{x\in S}\varphi_i(x),$

2. $|S^{\bot}|^{-\frac 12}\sum_{p\in
S^{\bot}}c_{\xi}(p){\varphi_{\xi}(\frac{p}{\xi})}=|S|^{-\frac
12}\sum_{x\in S}\varphi_{\xi}(x).$

\medskip
\noindent
When the generator of $S$ is $\varepsilon s$, we write this as the
following, explicitly :

\medskip

1. $H\exp(-i\frac{\pi}{4})\sum_{p\in
S^{\bot}}\exp(i\pi p^2)=s\sum_{x\in S}\exp(-i\pi x^2),$

2. $H\sum_{p\in
S^{\bot}}c_{\xi}(p){\exp(-\frac{1}{\xi}\pi p^2)}=s\sum_{x\in
S}\exp(-\xi\pi x^2).$

\medskip
\noindent
We obtain the following proposition :

\medskip

{\bf Proposition 4.4}

\noindent
(i)  If $s=H$, then the generator of $S$ is 1 and
$S=S^{\bot}=L\cap\,^{\ast}{\bf Z}$. Hence 

\medskip

1. $\exp(-i\frac{\pi}{4})\sum_{p\in
L\cap\,^{\ast}{\bf Z}}\exp(i\pi p^2)=\sum_{x\in L\cap\,^{\ast}{\bf
Z}}\exp(-i\pi x^2),$

2. $\sum_{p\in
L\cap\,^{\ast}{\bf Z}}c_{\xi}(p){\exp(-\frac{1}{\xi}\pi p^2)}=\sum_{x\in
L\cap\,^{\ast}{\bf Z}}\exp(-\xi\pi x^2).$

\medskip
\noindent
We put the standard part of the above, we obtain :

\medskip

1. $\exp(-i\frac{\pi}{4})\sum_{-\infty<p<\infty}\exp(i\pi
p^2)=\sum_{-\infty<x<\infty}\exp(-i\pi x^2),$

2. st$(\sum_{p\in
L\cap\,^{\ast}{\bf Z}}c_{\xi}(p){\exp(-\frac{1}{\xi}\pi
p^2)})=$st$(\sum_{x\in L\cap\,^{\ast}{\bf Z}}\exp(-\xi\pi
x^2))=\sum_{-\infty<n<\infty}\exp(-\xi\pi n^2)=\theta(i\xi)$, where
$\theta(z)$ is the $\theta$-function.

\medskip
\noindent
(ii)  If $\varepsilon s$ is infinitesimal, then the equation : 2.
$H\sum_{p\in S^{\bot}}c_{\xi}(p){\exp(-\frac{1}{\xi}\pi p^2)}=s\sum_{x\in
S}\exp(-\xi\pi x^2)$ implies the following :

\medskip 

2. st$(\sum_{p\in
S^{\bot}}c_{\xi}(p){\exp(-\frac{1}{\xi}\pi p^2)})=$st$(\varepsilon
s\sum_{x\in S}\exp(-\xi\pi x^2))$

$=\int_{-\infty}^{\infty}\exp(-\xi\pi
x^2)dx=\frac{1}{\sqrt{\xi}}.$

\medskip

It is known that
st($c_{\xi}(p)=\frac{1}{\sqrt{\xi}}$,
and $\sum_{-\infty<x<\infty}\exp(-\xi\pi x^2)$ in 2 of (i) is equal to
$\frac{1}{\sqrt{\xi}}\sum_{-\infty<p<\infty}{\exp(-\frac{1}{\xi}\pi
p^2)}$ by the standard Poisson summation formula. Hence, by 2 of (i),
st$(\sum_{p\in S^{\bot}}c_{\xi}(p){\exp(-\frac{1}{\xi}\pi
p^2)})=\sum_{-\infty<p<\infty}$st$(c_{\xi}(p){\exp(-\frac{1}{\xi}\pi
p^2)}).$

\medskip

We extend the above formulation of $\varphi_i(x)$ to
$\varphi_{im}(x)=\exp(-im\pi x^2)$, for an integer $m$ so that $m|2H^2$
. If $m|\frac {p}{\varepsilon}$, we recall

\medskip

$(F\varphi_{im})(p)=c_{im}(p)\exp(i\pi\frac 1m p^2)$, 

\medskip
\noindent
where
$c_{im}(p)=\sqrt{\frac m2}\frac{1+i^{\frac{2H^2}{m}}}{1+i}$ for a
positive $m$ and $c_{im}(p)=\sqrt{\frac
{-m}{2}}\frac{1+(-i)^{\frac{2H^2}{-m}}}{1-i}$ for a negative $m$.

Hence $|S^{\bot}|^{-\frac 12}\sum_{p\in
S^{\bot}}c_{im}(p){\varphi_{\frac{1}{im}}(p)}=|S|^{-\frac 12}\sum_{x\in S}
\varphi_{im}(x).$ When the generator $\varepsilon s'$ of $S^{\bot}$
satifies $m|s'$, that is, the generator $\varepsilon s$ of $S$ satifies
$m|\frac{H^2}{s}$, it reduces to the following :

\medskip

$H\sqrt{\frac m2}\frac{1+i^{\frac{2H^2}{m}}}{1+i}\sum_{p\in
S^{\bot}}\exp(i\pi\frac 1m p^2)=s\sum_{x\in
S}\exp(-im\pi x^2)$ for a positive $m$,  

$H\sqrt{\frac
{-m}{2}}\frac{1+(-i)^{\frac{2H^2}{-m}}}{1-i}\sum_{p\in
S^{\bot}}\exp(i\pi\frac 1m p^2)=s\sum_{x\in
S}\exp(-im\pi x^2)$ for a negative $m$.

\medskip
\noindent
If $s=H$ and $m|H$, then

\medskip

$\sqrt{\frac
m2}\frac{1+i^{\frac{2H^2}{m}}}{1+i}\sum_{-\infty<p<\infty}\exp(i\pi\frac
1m p^2)=\sum_{-\infty<x<\infty}\exp(-im\pi x^2)$

\medskip
\noindent
for a positive $m$,  

\medskip

$\sqrt{\frac
{-m}{2}}\frac{1+(-i)^{\frac{2H^2}{-m}}}{1-i}\sum_{-\infty<p<\infty}\exp(i\pi\frac
1m p^2)=\sum_{-\infty<x<\infty}\exp(-im\pi x^2)$ 

\medskip
\noindent
for a negative $m$, that is,

\medskip

$\sqrt{m}\exp(-i\frac{\pi}{4})\sum_{-\infty<p<\infty}\exp(i\pi\frac
1m p^2)=\sum_{-\infty<x<\infty}\exp(-im\pi x^2)$ for a positive $m$,  

$\sqrt{-m}\exp(i\frac{\pi}{4})\sum_{-\infty<p<\infty}\exp(i\pi\frac
1m p^2)=\sum_{-\infty<x<\infty}\exp(-im\pi x^2)$ for a negative $m$.

\medskip

We remark that it does not coincide with the formula
$\exp(i\frac{\pi}{4})\frac{1}{\sqrt{-z}}\theta(-\frac 1z)=\theta(z)$ for
$\theta$-function of Im$(z)>0$. The reason is that the above nonstandard
calculation implies an m multiple of the domain for the function
$\exp(-im\pi x^2)$.

\bigskip

{\bf 5. Poisson summation formula for
Definition 1.2 on the space of functionals}

\medskip

We extend Poisson summation formula of finite group to our 
infinitesimal Fourier transformation, Definition 1.2, on the space of
functionals originally defined in [N-O1]. 

\medskip

{\bf 5-1.  Formulation}

\medskip

{\bf Theorem 5.1.}\;\;Let $Y$ be an internal subgroup of $X$. Then we
obtain, for $f\in A$,

\medskip

$|Y^\bot|^{-\frac 12}\sum_{b\in Y^\bot}(Ff)(b)=|Y|^{-\frac 12}\sum_{a\in
Y}f(a),$

\medskip
\noindent
where $Y^\bot :=\{b\in X\,|\, \exp(2\pi i<a,b>)=1\textrm{ for }\forall
a\in X\}$ and $<a,b>:=\sum_{k\in L}a(k)b(k).$

\medskip

{\bf Lemma 5.2.}$\;\;|Y^\bot| =\frac{|X|}{|Y|}.$

\medskip

{\bf Proof of Lemma 5.2.}  For $k\in L$, we denote $Y_k:=\{a(k)\in L'\,|\,
a\in Y\}.$

\medskip

$b\in Y^\bot \Longleftrightarrow \forall a\in Y,\; \exp(2\pi i\sum_{k\in
L}a(k)b(k))=1$

$\Longleftrightarrow \forall a\in Y,\; \prod_{k\in L}(\exp(2\pi i
a(k)b(k)))=1$

$\Longleftrightarrow \forall k\in L,\; \forall a(k)\in Y_k,\; \exp(2\pi i
a(k)b(k))=1$

$\Longleftrightarrow \forall k\in L,\; b(k)\in Y_k^\bot$

$\Longleftrightarrow b:L\to L',\; \forall k\in L,\; b(k)\in Y_k^\bot.$

\medskip
\noindent
Hence $|Y^\bot|=\prod_{k\in L}|Y_k^\bot|.$ Theorem 3.1
implies $|Y_k^\bot|=\frac{H'{}^2}{|Y_k|}.$ Thus

\medskip

$|Y^\bot|=\prod_{k\in
L}\left(\frac{H'{}^2}{|Y_k|}\right)=\frac{H'{}^{2\,^{\ast}H}}{\prod_{k\in
L}|Y_k|}=\frac{|X|}{|Y|}.$

\medskip

{\bf Proof of Theorem 5.1.}

\medskip

$|Y^\bot|^{-\frac 12}\sum_{b\in Y^\bot}(Ff)(b)$

$=|Y^\bot|^{-\frac 12}\sum_{b\in Y^\bot}\sum_{a\in X}\varepsilon_0\exp(-2\pi
i<a,b>)f(a)$

$=|Y^\bot|^{-\frac 12}\sum_{a\in X}\varepsilon_0(\sum_{b\in
Y^\bot}\exp(-2\pi i<a,b>))f(a).$

\medskip
\noindent
$\textrm{Since }\sum_{b\in Y^\bot}\exp(-2\pi i<a,b>)=
\begin{cases}
0& \text{$(a\notin Y)$}\\
|Y^\bot|  & \text{$(a\in Y)$}
\end{cases}$, the above is equal to

\medskip
$|Y^\bot|^{-\frac 12}\varepsilon_0 |Y^\bot|\sum_{a\in Y}f(a)
=|Y^\bot|^{-\frac 12}(H')^{-\,^{\ast}H^2}\sum_{a\in Y}f(a)
=|Y|^{-\frac 12}\sum_{a\in Y}f(a). $

\medskip

Especially if $f$ is written as $\prod_{k\in L}f_k$, that is,
$f(a)=\prod_{k\in L}f_k(a(k))$, then $(Ff)(b)$ is $\sum_{a\in
X}\varepsilon_0\exp(-2\pi i\sum_{k\in L}a(k)b(k))\prod_{k\in
L}f_k(a(k))$. It is calculated to 

\noindent
$\prod_{k\in
L}(\sum_{a(k)\in
L'}\varepsilon'\exp(-2\pi ia(k)b(k))f_k(a(k))$, that represents an
infinite product of infinitesimal Fourier transformation defined by
Kinoshita. In general, since $f$ is not written as $\prod_{k\in
L}f_k$, our infinitesimal Fourier transformation is not represented as an
product of infinitesimal Fourier transformation defined by Kinoshita.

We summarize the argument, we obtain : $Ff=\prod_{k\in
L}F_k f_k$, where $F_k$ is an infinitesimal Fourier transformation for
each $k\in L$. We apply Proposition 5.3 to each $F_k$.

\medskip

{\bf Corollary 5.3}

\noindent
(i) If each generator of $Y_k$ is equal to $1$,
$f$ is written as $\prod_{k\in
L}f_k$, $f_k=\,^{\ast}($st$(f_k))|_{L'}$, and $\sum_{-\infty
<n<\infty}$st$(f_k)(n)$ converges, then 

\medskip

st$(\sum_{b\in
Y^\bot}(Ff)(b))=\prod_{k\in L}(\sum_{-\infty <n<\infty}$st$(f_k)(n))$.

\medskip
\noindent
(ii) If each generator of $Y_k$ is infinitesimal, $f$ is  written as $\prod_{k\in
L}f_k$, $f_k=\,^{\ast}($st$(f_k))|_{L'}$ and st$(f_k)$ is
$L_{1}$-integrable on  {\bf R}, then 

\medskip

st$(\sum_{b\in
Y^\bot}(Ff)(b))=\prod_{k\in L}\int_{-\infty <t<\infty}$st$(f_k)(t)dt$.

\bigskip

{\bf 5-2.  Examples}

\medskip

From now on the infinitesimal Fourier transformation $F({\bf
e}^{\sharp}(f))$ for a functional $f\in A$ is often denoted simply $Ff$.
We apply Theorem 5.2 to the following two functionals :

\medskip

1. $f_i(a)=\exp(-i\pi\sum_{k\in L}a(k)^2),$

2. $f_{\xi}(a)=\exp(-\xi\pi\sum_{k\in L}a(k)^2),$ where $\xi\in{\bf C}$,
Re$(\xi)>0.$

\medskip
\noindent
The infinitesimal Fourier transformations of the functionals are :

\medskip

1. $(Ff_i)(b)=(-1)^{\frac H2}\overline{f_i(b)}$,

2. $(Ff_{\xi})(b)=B_{\xi}(b){f_{\xi}(\frac{b}{\xi})}$,

\medskip
\noindent
hence we obtain :

\medskip

1. $|Y^{\bot}|^{-\frac 12}(-1)^{\frac H2}\sum_{b\in
Y^{\bot}}\overline{f_i(b)}=|Y|^{-\frac 12}\sum_{a\in Y}f_i(a),$

2. $|Y^{\bot}|^{-\frac 12}\sum_{b\in
Y^{\bot}}B_{\xi}(b){f_{\xi}(\frac{b}{\xi})}=|Y|^{-\frac
12}\sum_{a\in Y}f_{\xi}(a).$

\medskip
\noindent
We write this as the following, explicitly :

\medskip

1. $|Y^{\bot}|^{-\frac 12}(-1)^{\frac H2}\sum_{b\in
Y^{\bot}}\exp(-i\pi\sum_{k\in L}b(k)^2)=|Y|^{-\frac
12}\sum_{a\in Y}\exp(-i\pi\sum_{k\in L}a(k)^2),$

2. $|Y^{\bot}|^{-\frac 12}\sum_{b\in
Y^{\bot}}B_{\xi}(b)\,{\exp(-\frac{1}{\xi}\pi\sum_{k\in
L}b(k)^2)}=|Y|^{-\frac 12}\sum_{a\in
Y}\exp(-\xi\pi\sum_{k\in L}a(k)^2).$

\medskip

Corollaly 5.3 implies the following proposition 5.4.

\medskip

{\bf Proposition 5.4}

\noindent 
(i) If each generator of $Y_k$ is equal to $1$,
then

\medskip

1. $(-1)^{\frac H2}$st$(\sum_{b\in
Y^\bot}\exp(-i\pi\prod_{k\in
L}b(k)^2))=(\sum_{-\infty <n<\infty}\exp(-i\pi n^2))^{H^2}$,

2. st$(\sum_{b\in
Y^\bot}B_{\xi}(b)\exp(-\frac{1}{\xi}\pi\sum_{k\in
L}b(k)^2))=(\sum_{-\infty <n<\infty}\exp(-\xi\pi
n^2))^{H^2}$

$\left(=(\theta(i\xi))^{H^2}\right)$,

\medskip
\noindent 
(ii) If each generator of $Y_k$ is equal to a natural number $m_k$,
then

\medskip

1. $(-1)^{\frac H2}$st$(\sum_{b\in
Y^\bot}\exp(-i\pi\prod_{k\in
L}b(k)^2))=\prod_{k\in L}(m_k\sum_{-\infty <n<\infty}\exp(-i\pi
m_k^2 n^2))$,

2. st$(\sum_{b\in
Y^\bot}B_{\xi}(b)\exp(-\frac{1}{\xi}\pi\sum_{k\in
L}b(k)^2))=\prod_{k\in L}(m_k\sum_{-\infty <n<\infty}\exp(-\xi\pi m_k^2
n^2))$

$\left(=\prod_{k\in L}(m_k\theta(im_k^2 \xi))\right)$,

\medskip
\noindent 
(iii) If each generator of $Y_k$ is infinitesimal,
then

\medskip

2. st$(\sum_{b\in
Y^\bot}B_{\xi}(b)\exp(-\frac{1}{\xi}\pi\sum_{k\in
L}b(k)^2))=(\int_{-\infty}^{\infty}\exp(-\xi\pi
t^2)dt)^{H^2}$

$\left(=\left(\ast\left(\frac{1}{\sqrt{\xi}}\right)\right)^{H^2}\right)$.

\medskip

We extend the above formulation of $g_i(a)$ to
$g_{im}(a)=\exp(-im\pi \sum_{k\in L}a^2(k))$, for an integer $m$ so that
$m|2H'{}^2$ . If $m|\frac {b(k)}{\varepsilon'}$, we recall

$(Fg_{im})(b)=B_{im}(b)g_{\frac{1}{im}}(b)$, where
$B_{im}(b)=\left(\sqrt{\frac
m2}\frac{1+i^{\frac{2H'{}^2}{m}}}{1+i}\right)^{(\,^{\star}H)^2}$ for a positive
$m$ and $B_{im}(b)=\left(\sqrt{\frac
{-m}{2}}\frac{1+(-i)^{\frac{2H'{}^2}{-m}}}{1-i}\right)^{(\,^{\star}H)^2}$ for a
negative $m$.

Hence $|Y^{\bot}|^{-\frac 12}\sum_{b\in
Y^{\bot}}B_{im}(b){g_{\frac{1}{im}}(b)}=|Y|^{-\frac 12}\sum_{a\in Y}
g_{im}(a).$ When each generator $\varepsilon' s'{}_k$ of $Y^{\bot}_k$
satisfies $m|s'{}_k$, that is, each generator $\varepsilon' s{}_k$ of
$Y_k$ satisfies
$m|\frac{H'{}^2}{s{}_k}$, it reduces to the following :

\medskip

$H'{}^{(\,^{\star}H)^2}\left(\sqrt{\frac
m2}\frac{1+i^{\frac{2H'{}^2}{m}}}{1+i}\right)^{(\,^{\star}H)^2}\sum_{b\in
Y^{\bot}}
\exp(i\pi\frac 1m \sum_{k\in L}b(k)^2)$

$=\prod_{k\in L}s{}_k\sum_{a\in
Y}\exp(-im\pi \sum_{k\in L}a(k)^2)$ for a positive $m$, and 

$H'{}^{(\,^{\star}H)^2}\left(\sqrt{\frac
{-m}{2}}\frac{1+(-i)^{\frac{2H'{}^2}{-m}}}{1-i}\right)^{(\,^{\star}H)^2}
\sum_{b\in
Y^{\bot}}
\exp(i\pi\frac 1m \sum_{k\in L}b(k)^2)$

$=\prod_{k\in L}s{}_k\sum_{a\in
Y}\exp(-im\pi \sum_{k\in L}a(k)^2)$ for a negative $m$.

\medskip
\noindent
If $s{}_k=H'$ and $m|H'$, then

\medskip

$\left(\sqrt{\frac
m2}\frac{1+i^{\frac{2H'{}^2}{m}}}{1+i}\right)^{(\,^{\star}H)^2}\sum_{b\in
Y^{\bot}}
\exp(i\pi\frac 1m \sum_{k\in L}b(k)^2)$

$=\sum_{a\in
Y}\exp(-im\pi \sum_{k\in L}a(k)^2)$ for a positive $m$, and 

$\left(\sqrt{\frac
{-m}{2}}\frac{1+(-i)^{\frac{2H'{}^2}{-m}}}{1-i}\right)^{(\,^{\star}H)^2}
\sum_{b\in
Y^{\bot}}
\exp(i\pi\frac 1m \sum_{k\in L}b(k)^2)$

$=\sum_{a\in
Y}\exp(-im\pi \sum_{k\in L}a(k)^2)$ for a negative $m$, that is,

$\left(\sqrt{m}\exp(-i\frac{\pi}{4})\right)^{(\,^{\star}H)^2}\sum_{b\in
Y^{\bot}}
\exp(i\pi\frac 1m \sum_{k\in L}b(k)^2)$

$=\sum_{a\in
Y}\exp(-im\pi \sum_{k\in L}a(k)^2)$ for a positive $m$, and 

$\left(\sqrt{-m}\exp(i\frac{\pi}{4})\right)^{(\,^{\star}H)^2}\sum_{b\in
Y^{\bot}}
\exp(i\pi\frac 1m \sum_{k\in L}b(k)^2)$

$=\sum_{a\in
Y}\exp(-im\pi \sum_{k\in L}a(k)^2)$ for a negative $m$.

\bigskip

{\bf 6. Poisson summation formula for
Definition 1.3 on the space of functionals}

\medskip

We extend Poisson summation formula of finite group to our 
infinitesimal Fourier transformation, Definition 1.3, on the space of
functionals originally defined in [N-O1]. 

\medskip

{\bf 6-1. Formulation}

\medskip 

We obtain the following theorem for Definition 1.3 as the above argument.

\bigskip

{\bf Theorem 6.1.}\;\;Let $Y$ be an internal subgroup of $X$. Then we
obtain, for $f\in A$,

\medskip

$|Y^{\bot\varepsilon}|^{-\frac 12}\sum_{b\in
Y^{\bot\varepsilon}}(Ff)(b)=|Y|^{-\frac 12}\sum_{a\in Y}f(a),$

\medskip
\noindent
where $Y^{\bot\varepsilon} :=\{b\in X\,|\, \exp(2\pi
i<a,b>_{\varepsilon})=1\textrm{ for }\forall a\in X\}$
and $<a,b>_{\varepsilon} :=\,^{\star}\varepsilon \sum_{k\in L}a(k)b(k).$

\medskip

{\bf Lemma 6.2.}\;\;$|Y^{\bot\varepsilon}| =\frac{|X|}{|Y|}.$

\medskip

{\bf Proof of Lemma 6.2.}  For $k\in L$, we denote $Y_k:=\{a(k)\in L'\,|\,
a\in Y\}.$

\medskip

$b\in Y^{\bot\varepsilon} \Longleftrightarrow \forall a\in Y,\; \exp(2\pi
i\,^{\star}\varepsilon\sum_{k\in L}a(k)b(k))=1$

$\Longleftrightarrow \forall a\in Y,\; \prod_{k\in L}(\exp(2\pi
i\,^{\star}\varepsilon a(k)b(k)))=1$

$\Longleftrightarrow \forall k\in L,\; \forall a(k)\in Y_k,\; \exp(2\pi
i\,^{\star}\varepsilon a(k)b(k))=1$

$\Longleftrightarrow \forall k\in L,\; \,^{\star}\varepsilon b(k)\in
Y_k^\bot$

$\Longleftrightarrow b:L\to L',\; \forall k\in L,\; \,^{\star}\varepsilon
b(k)\in Y_k^\bot.$

\medskip
\noindent
For $k\in L$, we write $m$, $n$ as gereraters defined by :

\medskip

$Y_k=<\varepsilon' m>,\; \{b(k)\in L'\,|\, \,^{\star}\varepsilon b(k)\in
Y_k^\bot\}=<\varepsilon' n>.$

\medskip
\noindent
Now

$\exp(2\pi i\,^{\star}\varepsilon\varepsilon' m\varepsilon' n)=1
\Longleftrightarrow \,^{\star}\varepsilon\varepsilon' m\varepsilon' n
\in
\,^{\star}(\,^{\ast}{\bf Z})$

$\Longleftrightarrow\,^{\star}\varepsilon\varepsilon' m\varepsilon'
n=1\;\;\;\;\;\cdots \textrm{(1)}.$

\medskip
\noindent
We write $Y_k^{\bot\varepsilon}:=\{b(k)\in L'\,|\, \,^{\star}\varepsilon
b(k)\in Y_k^\bot\}.$ Then $|Y_k^{\bot\varepsilon}|=\frac{|L'|}{n}
=\frac{\,^{\star}HH'{}^2}{n}
=\frac{1}{\,^{\star}\varepsilon \varepsilon'{}^2n}
=m.
$ This is equal to $\frac{\,^{\star}HH'{}^2}{\frac{\,^{\star}HH'{}^2}{m}}
=\frac{|L'|}{|Y_k|}.$ 
Hence

\medskip

$|Y^{\bot\varepsilon}|=\prod_{k\in L}|Y_k^{\bot\varepsilon}|
=\prod_{k\in L}\frac{|L'|}{|Y_k|}
=\frac{|L'|^{\,^{\star}H^2}}{\prod_{k\in L}|Y_k|}
=\frac{|X|}{|Y|}.$

\medskip

{\bf Proof of Theorem 6.1.}

\medskip

$|Y^{\bot\varepsilon}|^{-\frac 12}\sum_{b\in Y^{\bot\varepsilon}}(Ff)(b)$

$=|Y^{\bot\varepsilon}|^{-\frac 12}\sum_{b\in Y^{\bot\varepsilon}}\sum_{a\in X}\varepsilon_0\exp(-2\pi
i<a,b>_\varepsilon)f(a)$

$=|Y^{\bot\varepsilon}|^{-\frac 12}\sum_{a\in X}\varepsilon_0(\sum_{b\in
Y^{\bot\varepsilon}}\exp(-2\pi i<a,b>_\varepsilon))f(a).$

\medskip
\noindent
$\textrm{Since }\sum_{b\in Y^{\bot\varepsilon}}\exp(-2\pi
i<a,b>_\varepsilon)=
\begin{cases}
0& \text{$(a\notin Y)$}\\
|Y^{\bot\varepsilon}|  & \text{$(a\in Y)$}
\end{cases}$, the above is equal to

\medskip

$|Y^{\bot\varepsilon}|^{-\frac 12}\varepsilon_0 |Y^{\bot\varepsilon}|\sum_{a\in Y}f(a)
=|Y^{\bot\varepsilon}|^{-\frac 12}(H')^{-\,^{\ast}H^2}\sum_{a\in Y}f(a)
=|Y|^{-\frac 12}\sum_{a\in Y}f(a).
$

\medskip

We obtain the following :

\medskip

{\bf Corollary 6.3}

\noindent
(i) If each generator of $Y_k$ is equal to $1$,
$f$ is written as $\prod_{k\in
L}f_k$, $f_k=\,^{\ast}($st$(f_k))|_{L'}$, and $\sum_{-\infty
<n<\infty}$st$(f_k)(n)$ converges, then

\medskip

$H^{\frac{H^2}{2}}$st$(\sum_{b\in Y^\bot}(Ff)(b))=\prod_{k\in
L}(\sum_{-\infty <n<\infty}$st$(f_k)(n))$.

\medskip
\noindent
(ii) If each generator of $Y_k$ is infinitesimal, $f$ is  written as $\prod_{k\in
L}f_k$, $f_k=\,^{\ast}($st$(f_k))|_{L'}$, and st$(f_k)$ is
$L_{1}$-integrable on  {\bf R}, then

\medskip

$H^{\frac{H^2}{2}}$st$(\sum_{b\in Y^\bot}(Ff)(b))=\prod_{k\in
L}\int_{-\infty}^{\infty}$st$(f_k)(t)dt$.

\bigskip

{\bf 6-2. Examples}

\medskip

We apply Theorem 3.3 to the following two functionals :

\medskip

1. $g_i(a)=\exp(-i\pi\,^{\star}\varepsilon\sum_{k\in L}a(k)^2),$

2. $g_{\xi}(a)=\exp(-\xi\pi\,^{\star}\varepsilon\sum_{k\in L}a(k)^2).$

\medskip
\noindent
The infinitesimal Fourier transformations of the functionals are :

\medskip

1. $(Fg_i)(b)=(-1)^{\frac H2}\overline{g_i(b)}$,

2. $(Fg_{\xi})(b)=C_{\xi}(b){g_{\xi}(\frac{b}{\xi})}$,

\medskip
\noindent
hence we obtain :

\medskip

1. $|Y^{\bot\varepsilon}|^{-\frac 12}(-1)^{\frac H2}\sum_{b\in
Y^{\bot\varepsilon}}\overline{g_i(b)}=|Y|^{-\frac 12}\sum_{a\in Y}g_i(a),$

2. $|Y^{\bot\varepsilon}|^{-\frac 12}\sum_{b\in
Y^{\bot\varepsilon}}C_{\xi}(b){g_{\xi}(\frac{b}{\xi})}=|Y|^{-\frac
12}\sum_{a\in Y}g_{\xi}(a).$

\medskip
\noindent
We write this as the following, explicitly :

\medskip

1. $|Y^{\bot\varepsilon}|^{-\frac 12}(-1)^{\frac H2}\sum_{b\in
Y^{\bot\varepsilon}}\exp(-i\pi\,^{\star}\varepsilon\sum_{k\in
L}b(k)^2)=|Y|^{-\frac 12}\sum_{a\in
Y}\exp(-i\pi\,^{\star}\varepsilon\sum_{k\in L}a(k)^2),$

2. $|Y^{\bot\varepsilon}|^{-\frac 12}\sum_{b\in
Y^{\bot\varepsilon}}C_{\xi}(b){\exp(-\frac{1}{\xi}\pi\,^{\star}\varepsilon\sum_{k\in
L}a(k)^2)}=|Y|^{-\frac 12}\sum_{a\in
Y}\exp(-\xi\pi\,^{\star}\varepsilon\sum_{k\in L}a(k)^2).$

\medskip

Corollaly 5.3 implies the following proposition 5.8.

\medskip

{\bf Proposition 6.4}

\noindent 
(i) If each generator of $Y_k$ is equal to $1$,
then the standard parts are :

\medskip

1. $H^{\frac{H^2}{2}}(-1)^{\frac H2}$st$(\sum_{b\in
Y^\bot_{\varepsilon}}\exp(-i\pi\varepsilon\sum_{k\in
L}b(k)^2))=(\sum_{-\infty <n<\infty}\exp(-i\pi\varepsilon n^2))^{H^2}$,

2. $H^{\frac{H^2}{2}}$st$(\sum_{b\in
Y^\bot_{\varepsilon}}C_{\xi}(b)\exp(-\frac{1}{\xi}\pi\varepsilon\sum_{k\in
L}b(k)^2))=(\sum_{-\infty <n<\infty}\exp(-\xi\pi\varepsilon
n^2))^{H^2}$

$\left(=(\theta(i\xi))^{H^2}\right)$,

\medskip
\noindent 
(ii) If each generator of $Y_k$ is equal to a natural number $m_k$,
then

\medskip

1. $H^{\frac{H^2}{2}}(-1)^{\frac H2}$st$(\sum_{b\in
Y^\bot_{\varepsilon}}\exp(-i\pi\varepsilon\sum_{k\in
L}b(k)^2))=\prod_{k\in L}(m_k\sum_{-\infty <n<\infty}\exp(-i\pi\varepsilon
m_k^2 n^2))$,

2. $H^{\frac{H^2}{2}}$st$(\sum_{b\in
Y^\bot_{\varepsilon}}C_{\xi}(b)\exp(-\frac{1}{\xi}\pi\varepsilon\sum_{k\in
L}b(k)^2))=\prod_{k\in L}(m_k\sum_{-\infty <n<\infty}\exp(-\xi\pi\varepsilon m_k^2
n^2))$

$\left(=\prod_{k\in L}(m_k\theta(im_k^2 \xi))\right)$,

\medskip
\noindent 
(iii) If each generator of $Y_k$ is infinitesimal,
then

\medskip

2. st$(\sum_{b\in
Y^\bot_{\varepsilon}}C_{\xi}(b)\exp(-\frac{1}{\xi}\pi\varepsilon\sum_{k\in
L}b(k)^2))=(\int_{-\infty}^{\infty}\exp(-\xi\pi
t^2)dt)^{H^2}$

$\left(=\left(\ast\left(\frac{1}{\sqrt{\xi}}\right)\right)^{H^2}\right)$.
 
\medskip

We extend the above formulation of $g_i(a)$ to
$g_{im}(a)=\exp(-im\pi \,^{\star}\varepsilon\sum_{k\in L}a^2(k))$, for an
integer
$m$ so that
$m|2\,^{\star}HH'{}^2$ . If $m|\frac {b(k)}{\varepsilon'}$ for an
arbitrary $k\in L$, we recall

$(Fg_{im})(b)=C_{im}(b)g_{\frac{1}{im}}(b)$, where
$C_{im}(b)=\left(\sqrt{\frac
m2}\frac{1+i^{\frac{2\,^{\star}HH'{}^2}{m}}}{1+i}\right)^{\,^{\star}H^2}$
for a positive $m$ and $C_{im}(b)=\left(\sqrt{\frac
{-m}{2}}\frac{1+(-i)^{\frac{2\,^{\star}HH'{}^2}{-m}}}{1-i}\right)^{\,^{\star}H^2}$
for a negative $m$.

Hence $|Y^{\bot_{\varepsilon}}|^{-\frac 12}\sum_{b\in
Y^{\bot}}C_{im}(b){g_{\frac{1}{im}}(b)}=|Y|^{-\frac 12}\sum_{a\in Y}
g_{im}(a).$ When each generator $\varepsilon' s'{}_k$ of
$Y^{\bot_{\varepsilon}}_k$ satisfies $m|s'{}_k$, that is, each generator
$\varepsilon' s{}_k$ of
$Y_k$ satisfies
$m|\frac{\,^{\star}HH'{}^2}{s{}_k}$, it reduces to the following :

\medskip

$H^{\frac{H^2}{2}}H'{}^{(\,^{\star}H)^2}\left(\sqrt{\frac
m2}\frac{1+i^{\frac{2\,^{\star}HH'{}^2}{m}}}{1+i}\right)^{(\,^{\star}H)^2}\sum_{b\in
Y^{\bot_{\varepsilon}}}
\exp(i\pi\frac 1m \,^{\star}\varepsilon\sum_{k\in L}b(k)^2)$

$=\prod_{k\in
L}s{}_k\sum_{a\in Y}\exp(-im\pi \,^{\star}\varepsilon\sum_{k\in
L}a(k)^2)$ for a positive $m$, and 

$H^{\frac{H^2}{2}}H'{}^{(\,^{\star}H)^2}\left(\sqrt{\frac
{-m}{2}}\frac{1+(-i)^{\frac{2\,^{\star}HH'{}^2}{-m}}}{1-i}\right)^{
(\,^{\star}H)^2}\sum_{b\in
Y^{\bot_{\varepsilon}}}
\exp(i\pi\frac 1m \,^{\star}\varepsilon\sum_{k\in L}b(k)^2)$

$=\prod_{k\in L}s{}_k\sum_{a\in
Y}\exp(-im\pi \,^{\star}\varepsilon\sum_{k\in L}a(k)^2)$ for a negative $m$.

\medskip
\noindent
If $s{}_k=H'$ and $m|H'$, then

\medskip

$H^{\frac{H^2}{2}}\left(\sqrt{\frac
m2}\frac{1+i^{\frac{2\,^{\star}HH'{}^2}{m}}}{1+i}\right)^
{(\,^{\star}H)^2}\sum_{b\in
Y^{\bot_{\varepsilon}}}
\exp(i\pi\frac 1m \sum_{k\in L}b(k)^2)$

$=\sum_{a\in
Y}\exp(-im\pi \,^{\star}\varepsilon\sum_{k\in L}a(k)^2)$ for a positive $m$, and 

$H^{\frac{H^2}{2}}\left(\sqrt{\frac
{-m}{2}}\frac{1+(-i)^{\frac{2\,^{\star}HH'{}^2}{-m}}}{1-i}\right)^{(
\,^{\star}H)^2}\sum_{b\in
Y^{\bot_{\varepsilon}}}
\exp(i\pi\frac 1m \sum_{k\in L}b(k)^2)$

$=\sum_{a\in
Y}\exp(-im\pi \,^{\star}\varepsilon\sum_{k\in L}a(k)^2)$ for a negative $m$, that is,

\noindent
$H^{\frac{H^2}{2}}\left(\sqrt{m}\exp(-i\frac{\pi}{4})\right)^{(\,^{\star}H)^2}\sum_{b\in
Y^{\bot_{\varepsilon}}}
\exp(i\pi\frac 1m \sum_{k\in L}b(k)^2)$

$=\sum_{a\in
Y}\exp(-im\pi \,^{\star}\varepsilon\sum_{k\in L}a(k)^2)$ for a positive $m$, and 

$H^{\frac{H^2}{2}}\left(\sqrt{-m}\exp(i\frac{\pi}{4})\right)^{(\,^{\star}H)^2}\sum_{b\in
Y^{\bot_{\varepsilon}}}
\exp(i\pi\frac 1m \sum_{k\in L}b(k)^2)$

$=\sum_{a\in
Y}\exp(-im\pi \,^{\star}\varepsilon\sum_{k\in L}a(k)^2)$ for a negative $m$.

\bigskip

{\bf 7. The infinitesimal Fourier transformation of a
functional $Z_s(a)$}

\medskip

In this section, we define a functional on $X$, and study a relationship between
the functional and the Riemann zeta function. We order all prime numbers
as
$p(1)=2$,
$p(2)=3$, ... ,
$p(n)<p(n+1)$, ... , that is, $p$ is a mapping from ${\bf N}$ to the set
$\{$prime number$\}$, $p:{\bf N}\to\{$prime number$\}$. The nonstandard
extension $\,^{\ast}p:\,^{\ast}{\bf N}\to\,^{\ast}\{$prime number$\}$ is
written as $\,^{\ast}p([l_{\mu}])=[p(l_{\mu})]$, and we define a mapping
$\tilde{p}:\,^{\ast}{\bf N}\to \,^{\star}(\,^{\ast}\{$prime number$\})$
as $\tilde{p}([l_{\mu}])= \,^{\star}[p(l_{\mu})]$. For $s\in {\bf C}$, we
define $Z_s(\in A)$ as the following :

\medskip

$Z_s(a):=\prod_{k\in
L}\tilde{p}(H(k+\frac{H}{2})+1)^{(-s(a(k)+\frac{H'}{2}))}$,

\medskip
\noindent
now $H(k+\frac{H}{2})+1$ is an element of $\,^{\ast}{\bf N}$ and
$a(k)+\frac{H'}{2}$ is an element of
$\,^{\star}(\,^{\ast}{\bf N})$. Then $Z_s(a)$ is calculated as
$\exp(-s\sum_{k\in L}\log(
\tilde{p}(H(k+\frac{H}{2})+1))a(k))\prod_{k\in
L}\tilde{p}(H(k+\frac{H}{2})+1)^{-s\frac{H'}{2}}$. We obtain
the following theorem for the Fourier transformation of ${\bf
e}^{\sharp}(Z_s)$ for Definition 1.2 :

\medskip

{\bf Theorem 7.1.}$\;\;(F({\bf e}^{\sharp}(Z_s)))(b)
=\left(\prod_{k\in
L}\tilde{p}(H(k+\frac{H}{2})+1)\right)^{-s\frac{H'}{2}}\\
\cdot\prod_{k\in
L}\varepsilon'\frac{\sinh((2\pi
i
b(k)+s\log\tilde{p}(H(k+\frac{H}{2})+1))\frac{H'}{2})}
{\exp(-\frac {\varepsilon'}{2}(2\pi i
b(k)+s\log\tilde{p}(H(k+\frac{H}{2})+1))\sinh(\frac {\varepsilon'}{2}(2\pi i
b(k)+s\log\tilde{p}(H(k+\frac{H}{2})+1))}$.

\medskip

{\bf Proof.}$\;\;(F({\bf e}^{\sharp}(Z_s)))(b)=\left(\prod_{k\in
L}\tilde{p}(H(k+\frac{H}{2})+1)\right)^{-s\frac{H'}{2}}\\
\cdot\sum_{a\in X}\varepsilon_0\exp(-s\sum_{k\in
L}\log\tilde{p}(H(k+\frac{H}{2})+1)a(k))\exp(-2\pi i\,\sum_{k\in
L}a(k)b(k))$

$=\left(\prod_{k\in
L}\tilde{p}(H(k+\frac{H}{2})+1)\right)^{-s\frac{H'}{2}}\\
\cdot\sum_{a\in X}\varepsilon_0\exp(-(2\pi i\,
b(k)+s\log\tilde{p}(H(k+\frac{H}{2})+1))a(k))$

$=\left(\prod_{k\in
L}\tilde{p}(H(k+\frac{H}{2})+1)\right)^{-s\frac{\,^{\star}HH'{}^2}{2}}\\
\cdot\prod_{k\in L}\varepsilon'\sum_{a(k)\in
L'}\exp(-(2\pi i\,b(k)+s\log\tilde{p}(H(k+\frac{H}{2})+1))a(k))$

$=\left(\prod_{k\in
L}\tilde{p}(H(k+\frac{H}{2})+1)\right)^{-s\frac{H'}{2}}\\
\cdot\prod_{k\in
L}\varepsilon'\frac{\exp((2\pi
i\,
b(k)+s\log\tilde{p}(H(k+\frac{H}{2})+1))\frac{H'}{2})-
\exp(-(2\pi i\,
b(k)+s\log\tilde{p}(H(k+\frac{H}{2})+1))\frac{H'}{2})}
{1-\exp(-\varepsilon'(2\pi i\,^{\star}\varepsilon
b(k)+s\log\tilde{p}(H(k+\frac{H}{2})+1))}$

$=\left(\prod_{k\in
L}\tilde{p}(H(k+\frac{H}{2})+1)\right)^{-s\frac{H'}{2}}\\
\cdot\prod_{k\in
L}\varepsilon'
\frac{\sinh((2\pi
i\,
b(k)+s\log\tilde{p}(H(k+\frac{H}{2})+1))\frac{H'}{2})}
{\exp(-\frac {\varepsilon'}{2}(2\pi i\,
b(k)+s\log\tilde{p}(H(k+\frac{H}{2})+1))\sinh(\frac {\varepsilon'}{2}(2\pi
i\, b(k)+s\log\tilde{p}(H(k+\frac{H}{2})+1))}.$

\medskip

We denote the Riemann zeta function by $\zeta(s)$, defined by
$\zeta(s)=\prod_{l=1}^{\infty}\frac{1}{1-p(l)^{-s}}$ for Re$(s)>1$. Let $Y_{\bf
Z}$ be a subgroup of $X$ so that each generator of $(Y_{\bf
Z})_k$ is equal to 1. Then we obtain the following
theorem :

\medskip

{\bf Theorem 7.2.}$\;\;$If Re$(s)>1$, then
st(st($\sum_{a\in Y_{\bf
Z}}{\bf e}^{\sharp}(Z_s))(a)))=\zeta(s).$

\medskip

{\bf Proof.}$\;\;$
st(st($\sum_{a\in Y_{\bf
Z}}{\bf e}^{\sharp}(Z_s))(a)))$

$=$st$\Bigl($st$\Bigl(\left(\prod_{k\in
L}\tilde{p}(H(k+\frac{H}{2})+1)\right)^{(-s(a(k)+\frac{H'}{2}))}\Bigr)\Bigr)$

$=$
st$\Bigl($st$\Bigl(\prod_{k\in
L}\frac{1-
\tilde{p}(H(k+\frac{H}{2})+1){}^{-sH'}}
{1-\tilde{p}(H(k+\frac{H}{2})+1){}^{-s}}\Bigr)\Bigr)$

$\Bigl($st$\Bigl(\prod_{k\in
L}\frac{1}
{1-\tilde{p}(H(k+\frac{H}{2})+1){}^{-s}}\Bigr)\Bigr)=\zeta(s).$

\medskip

Furthermore, Corollary 5.3.(1) and Theorem 7.2 imply the following :
\medskip

{\bf Corollary 7.3.}$\;\;$st$(\sum_{b\in
Y^\bot_{\bf Z}}(F({\bf e}^{\sharp}(Z_s)))(b))=$
st$\Bigl(\prod_{k\in
L}\frac{1-
\tilde{p}(H(k+\frac{H}{2})+1){}^{-sH'}}
{1-\tilde{p}(H(k+\frac{H}{2})+1){}^{-s}}\Bigr).$

\medskip

Hence we obtain : st(st$(\sum_{b\in
Y^\bot_{\bf Z}}(F({\bf e}^{\sharp}(Z_s)))(b))))=\zeta(s)$ for Re$(s)>1$.

\bigskip

\textsc{Acknowledgement}.  We would like to thank Prof. R.
Kobayashi for a useful suggestion about Poisson summation formula.

\bigskip 

\begin{center}
{\bf References}
\end{center}

[F-H] \; R.P. Feynman, A.R. Hibbs, Quantum mechanics and path integrals,
McGrow-Hill Inc. All rights (1965).

[G] \; E.I. Gordon, Nonstandard methods in commutative harmonic analysis,
Translations of mathematical monographs {\bf 164} American mathematical
society, 1997.

[K1]	\; M. Kinoshita, Nonstandard representation of distribution I, Osaka J.
Math. {\bf 25} (1988), 805-824.

[K2] 	\; M. Kinoshita,: Nonstandard representation of distribution II. 
Osaka J. Math. {\bf 27} (1990), 843-861.

[N-O1]  \; T. Nitta and T. Okada, Double infinitesimal Fourier
transformation for the space of functionals and reformulation of Feynman path
integral, Lecture Note Series in Mathematics, Osaka University Vol.{\bf 7} (2002),
255-298 in
Japanese.

[N-O2]\; T. Nitta, T. Okada, Infinitesimal Fourier
transformation for the space of functionals, preprint.

[N-O-T]\; T. Nitta, T. Okada and A. Tzouvaras, Classification of
non-well-founded sets and an application, Math. Log. Quart. {\bf 49}
(2003), 187-200.

[R] \; R. Remmert, Theory of complex functions, Graduate Texts in Mathematics
{\bf 122}, Springer, Berlin-Heidelberg-New York, 1992. 

[Sai]\; M. Saito, Ultraproduct and non-standard analysis, in Japanese,
Tokyo tosho, 1976.

[Sat]  \; I. Satake, The temptation to algebra, in
Japanese, Yuseisha, 1996.

[T]\; G. Takeuti, Dirac space, Proc. Japan Acad. {\bf 38}
(1962), 414-418.

\bigskip
\bigskip

{\small Takashi NITTA

Department of Education

Mie University 

Kamihama, Tsu, 514-8507, Japan

e-mail : nitta@edu.mie-u.ac.jp

\bigskip

Tomoko OKADA

Graduate school of Mathematics

Nagoya University

Chikusa-ku, Nagoya, 464-8602, Japan

e-mail : m98122c@math.nagoya-u.ac.jp}

\end{document}